\documentclass{article}
\usepackage{amssymb}

\begin{document}

\title{Hodge-Frobenius equations and the Hodge-B\"acklund transformation}
\author{Antonella Marini\thanks{%
Permanent address:  Dipartimento di Matematica, Universit\`{a} di
L'Aquila, 67100 L'Aquila, Italy; email: marini@dm.univaq.it} $\,$
and
Thomas H. Otway\thanks{%
email: otway@yu.edu} \\
\\
\textit{Department of Mathematics, Yeshiva University,}\\
\textit{\ \ New York, New York 10033}}
\date{}
\maketitle

\begin{abstract}

Linear and nonlinear Hodge-like systems for 1-forms are studied,
with an assumption equivalent to complete integrability substituted
for the requirement of closure under exterior differentiation. The
systems are placed in a variational context and properties of
critical points investigated. Certain standard choices of energy
density are related by B\"acklund transformations which employ basic
properties of the Hodge involution. These \emph{Hodge-B\"acklund
transformations} yield invariant forms of classical B\"acklund
transformations that arise in diverse contexts. Some extensions to
higher-degree forms are indicated. \emph{MSC2000}: 35A15, 58J72,
58A14.

\medskip

\noindent\emph{Key words}: nonlinear Hodge theory, B\"acklund
transformation, Frobenius Theorem

\end{abstract}

\section{Introduction}

The study of vectors which are both divergence-free and curl-free
can be traced back at least to Helmholtz's analysis of vortices and
gradients \cite{H}. The generalization to differential forms which
are both closed and co-closed under exterior differentiation is the
content of the Hodge equations; see, \emph{e.g.}, \cite{M}, Ch.\ 7.
The divergence-free condition is frequently relaxed in variational
contexts, but generalizations of the curl-free condition remain
rather rare.

Our goal is to study both linear and nonlinear variants of the Hodge
equations for differential forms which are neither co-closed nor
closed, but which satisfy milder conditions having physical and
geometric significance.

\subsection{Organization of the paper}

Sections 1-4 are mainly expository. We introduce the topic in Sec.\
1.2 with an example from fluid dynamics. Section 2 presents the
equations in an invariant context. The linear case is studied in
Sec.\ 3, largely as motivation for the considerably more complex
nonlinear case. Two geometric analogies are discussed in Sec.\ 4.
Technical results on the properties of solutions are presented in
Sec.\ 5. Section 6 shows that ideas introduced, in very different
contexts, by Yang \cite{Y} and by Magnanini and Talenti
\cite{MT1}-\cite{MT3}, can be given a unified interpretation in
terms of equations studied in the preceding sections.

The proofs of Theorems 5, 6, and 8, which are based on rather
straightforward applications of nonlinear elliptic theory, are
collected in an appendix, Sec.\ 7. We note that these applications
are only straightforward once solutions have been associated to a
uniformly sub-elliptic operator; this is accomplished in Sec.\ 5.2.
The methods used in that section to derive uniform estimates are
elementary and, in particular, do not require a delicate limiting
argument that has become known as \emph{Shiffman regularization}
(\emph{c.f.} \cite{Sh} and the appendix by R. J. Sibner to
\cite{Si}).

\subsection{A motivating example: steady, ideal flow}

In models for the steady, adiabatic and isentropic flow of an ideal
fluid, conservation of mass is represented by the \emph{continuity
equation}
\begin{equation}\label{conti}
    \nabla\cdot\left(\rho\left(|\mathbf{v}|^2\right)\mathbf{v}\right)=0,
\end{equation}
where $\mathbf{v}$ denotes flow velocity and $\rho$ denotes mass
density. The dependence of $\rho$ on $|\mathbf{v}|^2$ is a
consequence of compressibility; in the incompressible limit, eq.\
(\ref{conti}) says merely that the vector $\mathbf{v}$ has zero
divergence. If the fluid is irrotational, then the velocity is
curl-free in the sense that
\begin{equation}\label{cf}
    \nabla\times\mathbf{v}=0.
\end{equation}
Condition (\ref{cf}) implies, by the Poincar\'e Lemma, that there
exists locally a scalar flow potential
$\varphi\left(\mathbf{x}\right),$ where $\mathbf{x}\in\mathbb{R}^3$
denotes the position of a particle in the flow.

Perhaps the mildest weakening of the irrotationality condition
results from replacing (\ref{cf}) by the integrability condition
\begin{equation}\label{inte}
    \mathbf{v}\cdot\nabla\times\mathbf{v}=0.
\end{equation}

The replacement of the linear condition (\ref{cf}) by the nonlinear
condition (\ref{inte}) as a side condition to eq.\ (\ref{conti}) is
likely to result in singular solutions, even in the subsonic regime.
The usual arguments for reducing $\rho$ to the conventional form,
which depend on smoothness (see, \emph{e.g.}, Ch.\ 1.2 of
\cite{Be2}), would not necessarily apply in such cases. This
suggests that we consider whether a useful \emph{a priori} bound can
be placed on the size of the singular set for solutions of systems
having the general form (\ref{conti}), (\ref{inte}). In Sec.\ 5.2 we
take the first step toward an answer to this difficult question,
deriving sufficient conditions under which a solution remains
bounded on an apparent singular set of given codimension. Despite
the physical motivation (here and in various other examples
scattered throughout the text), our main interest in this paper is
in deriving hypotheses which are mathematically natural and apply to
large classes of mass densities.

\section{An invariant formulation}

Thus we generalize the mathematical context of Sec.\ 1.2. Let
$\Omega$ be an open, finite domain of $\mathbb{R}^n,$ $n \geq 2,$
satisfying an interior sphere condition. Consider the system
(\cite{O1}, Sec.\ VI; \cite{O2}, Sec.\ 4)
\begin{equation}\label{HF1}
    \delta\left(\rho(Q)\omega\right)=0,
\end{equation}

\begin{equation}\label{HF2}
    d\omega = \Gamma\wedge\omega,
\end{equation}
for scalar-valued 1-forms $\omega$ and $\Gamma,$ where $\Gamma$ is
given and $\omega$ is unknown; $d$ is the flat exterior derivative
with formal adjoint $\delta;$ $Q$ is a quadratic form in $\omega$
given by
\begin{equation}\label{speed}
    Q(\omega)=*\left(\omega \wedge *\omega\right)\equiv\langle\omega,\omega\rangle,
\end{equation}
where $*\colon \Lambda^k\rightarrow\Lambda^{n-k}$ is the Hodge
involution; $\rho$ is a positive, continuously differentiable
function of $Q$ (but a possibly singular function of $\mathbf{x}$).

Using (\ref{HF2}), we find that
\begin{equation}\label{froblink}
\omega\wedge d\omega=-\Gamma\wedge\omega\wedge\omega=0.
\end{equation}
A 1-form $\omega$ that is the pointwise Riemannian inner product
with a vector field $\mathbf{v}$ is said to be \emph{dual} to
$\mathbf{v}.$ In this case the left-hand side of (\ref{froblink}) is
equivalent to the left-hand side of (\ref{inte}), and any solution
of $\omega$ of eq.\ (\ref{HF2}) is dual to a solution $\mathbf{v}$
of eq.\ (\ref{inte}).

The system (\ref{HF1}), (\ref{HF2}) is uniformly elliptic provided
the differential inequality
\begin{equation}\label{unison}
    0<\kappa_1\leq\frac{\left(d/dQ\right)\left[  Q\rho^{2}(Q)\right]
}{\rho\left( Q\right) }\leq\kappa_2<\infty
\end{equation}
is satisfied for constants $\kappa_1,\kappa_2.$ In the context of
fluid dynamics, one typically encounters the weaker condition
\begin{equation}\label{subsonic}
    0<\rho^2(Q)+2Q\rho'(Q)\rho(Q).
\end{equation}
Ideal flow governed by eq.\ (\ref{conti}) is subsonic provided
(\ref{subsonic}) is satisfied. Moreover, there is typically a
critical value $Q_{crit}$ such that the right-hand side of
(\ref{subsonic}) tends to zero in the limit as $Q$ tends to
$Q_{crit}.$ In this case, eqs.\ (\ref{HF1}), (\ref{HF2}) with
$\rho$ satisfying (\ref{subsonic}) are elliptic, but not uniformly
so, and this condition has mathematical as well as physical
interest.

If $\Gamma\equiv 0,$ then the system (\ref{HF1}), (\ref{HF2})
degenerates to the \emph{nonlinear Hodge equations} introduced in
\cite{SS1} on the basis of a conjecture in \cite{Be1}. In that case
condition (\ref{HF2}) generates a cohomology class, which is not
true in the more general case studied here.

We show in Sec.\ 2.1 that whenever $\omega$ is a 1-form, eq.\
(\ref{HF2}) possesses solutions of the form
\begin{equation}\label{solu}
    \omega = e^\eta du,
\end{equation}
where $u$ and $\eta$ are 0-forms, and $\Gamma$ can be made
\emph{exact}: we can write $\Gamma=d\eta.$ When $\omega$ is a
$k$-form, a more general representation applies; that representation
is discussed in Sec.\ 2.2.

In Sec.\ 5.1 we introduce a variant of eq.\ (\ref{HF1}) for
differential $k$-forms ($k\geq 1$)  satisfying eq.\ (\ref{HF2})
with $\Gamma$ exact ($\Gamma\equiv d\eta$), namely
\begin{equation}\label{HB}
    \delta\left[\rho(Q)\omega\right]=(-1)^{n(k+1)}\ast \left( d\eta \wedge \ast \rho(Q)\omega\right),
\end{equation}
which arises as a variational equation of the \emph{nonlinear Hodge
energy}
\begin{equation}\label{energy}
    E = \frac{1}{2}\int_\Omega\int_0^Q\rho(s)ds\,d\Omega .
\end{equation}

\subsection{The Frobenius Theorem}

Let $\Gamma$ be fixed, and define
\[
S\equiv\{ \omega\in \Lambda (\Omega)\equiv \bigoplus_{k=1}^n
\Lambda^k (\Omega)\;:\; d\omega=\Gamma\wedge \omega\}.
\]
Denote by $I\equiv I(S)$ the ideal generated by $S$. If $\Gamma\neq
0,$ then clearly $dI\neq \{0\},$ so the $k$-forms $\omega\in S$ do
not generate cohomology classes.

Nevertheless, the ideal $I$ is \emph{closed}, \emph{i.e.},
$dI\subset I$. In fact, a differential form  $\alpha\in I$ is a
linear combination of forms of type $\omega\wedge\beta$ with
$\omega\in S,$ $\beta\in \Lambda (\Omega)$. The latter satisfy
\[
d(\omega\wedge\beta)= \pm \omega\wedge(\Gamma\wedge\beta \pm
d\beta),
\]
and thus satisfy $d\alpha\in I$. This is an important fact,
especially for exterior systems of $1$-forms.

Following the approach in Sec.\ 4-2 of \cite{E}, we define an
exterior system $\lbrace \omega^a \rbrace,$ $a = 1, \ldots, r$ of
$r$ 1-forms in a space of dimension $n = r + s$ to be
\emph{completely integrable} if and only if there exist $r$
independent functions $g^a,$ $a= 1, \ldots, r$ such that each of the
1-forms $\omega^a$ vanishes on the $r$-parameter family of
$s$-dimensional hypersurfaces $\lbrace g^a=k^a\,,\; a=1,\ldots r
\rbrace $ generated by letting the constants $k^a$ range over all
$r$-tuples of real numbers.

Equivalently, we define $\lbrace \omega^a \rbrace_1^r$  to be
\emph{completely integrable} if and only if there exists a
nonsingular $r\times r$ matrix of functions $\xi^a_b,$ and $r$
independent functions ${g^b}$ such that
\[
\omega^a = \sum_{b=1}^r \xi^a_b \, dg^b\;.
\]

The \emph{Frobenius Theorem} asserts that an exterior system
$\lbrace \omega^a \rbrace_1^r$ of 1-forms is completely integrable
if and only if it generates a closed ideal of  $\Lambda (\Omega ).$

Because a 1-form $\omega$ satisfying (\ref{HF2}) generates a closed
ideal, by the Frobenius Theorem it can always be written in the form
(\ref{solu}). (In this case, $r=1.$) Thus
\begin{equation}\label{grad_recurs}
    d\omega = d\eta\wedge \omega,
\end{equation}
which shows that $\Gamma$ can be chosen to be exact. Notice the
gauge invariance having the form $\Gamma\to \tilde \Gamma\equiv
\Gamma + f(\mathbf{x})\,\omega.$

For this reason, we call the system (\ref{HF1}), (\ref{HF2}) the
\emph{nonlinear Hodge-Frobenius equations for $1$-forms.}

Unfortunately, the Frobenius Theorem does not generalize to forms of
arbitrary degree $k,$ as the condition $dI\subset I$ does not imply
complete integrability if $k\neq 1$. However, this does not mean
that there is nothing to be said about higher-degree forms. Relevant
properties of such forms are described in the following section.

\subsection{Recursive forms}

An exterior differential form of degree $k$ is said to be
\emph{recursive with coefficient $\Gamma$} if it satisfies
(\ref{HF2}). Let $\Omega$ be star-shaped. It is known that one can
define a \emph{homotopy operator} $\mathcal{H}\,:
\Lambda^k(\Omega)\to\Lambda^{k-1}(\Omega)$, which satisfies
\begin{equation} \label{H1}
    \omega = d\mathcal{H} \omega + \mathcal{H} d\omega.
\end{equation}
This property can be used, among other things, to show that a closed
form on a star-shaped domain is exact. We omit the formal definition
of this operator (for this and further details see Sec.\ 5-3 of
\cite{E}), and only describe its main properties:
\begin{itemize}
\item[(a)] $\mathcal{H}$ is linear \item[(b)] $\mathcal{H}^2=0$
\item[(c)] $\mathcal{H}d\mathcal{H}=\mathcal{H}$,
$d\mathcal{H}d=d$ \item[(d)] $(d\mathcal{H})^2 = d\mathcal{H}$,
$(\mathcal{H}d)^2= \mathcal{H}d.$
\end{itemize}
Using (b), we observe that $\mathcal{H}d\omega\in Ker \,
\mathcal{H}.$ This and (\ref{H1}) can be used to define the
\emph{exact part of $\omega$} as $\omega_e \equiv d\mathcal{H}
\omega$ and the \emph{anti-exact part of $\omega$} as
$\omega_a\equiv \omega - \omega_e = \mathcal{H} d\omega.$

Using (\ref{H1}), one can further show that $\mathcal{H}$ improves
regularity. With no loss of generality, we prove this for forms
having vanishing anti-exact part, \emph{i.e.,} forms $\omega$ such
that
\begin{equation} \label{exact}
    \omega =\omega_e.
\end {equation}
Note that for any given form $\omega,$ no cancelations can occur
between $\omega_e$ and $\omega_a,$ therefore $\omega_e$ is always as
smooth as $\omega$, and furthermore $\mathcal{H}\omega_a =0$. Thus
we can restrict our attention to forms satisfying (\ref{exact}). If
$\omega$ is a one-form, then $\mathcal{H}\omega$ is a function and
(\ref{exact}) implies that
\[
\omega_i = \frac{\partial (\mathcal{H}\omega)}{\partial
x^i}\,\,\forall i,
\]
thus improving regularity. In the general case,
$\mathcal{H}\omega\equiv (\mathcal{H}\omega)_I dx^I,$ where $I$ is a
multi-index satisfying $|I|= k-1$. In order to improve regularity in
this case one would need to control $\delta (\mathcal{H}\omega),$
\emph{i.e.,} all derivatives of type $\partial
(\mathcal{H}\omega)_I/\partial x^i$ for $i\in I.$ But $\delta
(\mathcal{H}\omega)=0$ from the Hodge Decomposition Theorem
\cite{M}, as $\mathcal{H}\omega$ is anti-exact (using $(b),$ above);
so $\mathcal{H}$ is a smoothing operator on $k$-forms.

The important result for us is the following:

Recursive $k$-forms with coefficient $\Gamma$ on a star-shaped
region can be represented as follows \cite{E}:
\begin{equation} \label{recursive}
    \omega = e^\eta\left[d u + \mathcal{H}(\theta \wedge du)\right]\;,
\end{equation}
where $\eta = \mathcal{H} \Gamma$, $\theta = \mathcal{H} d\Gamma$,
and $u = \mathcal{H} (e^{-\eta}\omega)$. Condition (\ref{HF2})
implies that $\theta$ satisfies
\[
d\theta\wedge\left[d u + \mathcal{H}(\theta \wedge du)\right]=0.
\]
For our purposes we can rewrite (\ref{recursive}) as
\begin{equation}\label{higher1}
    \omega = e^\eta g(d u),
\end{equation}
where $g$ is a smooth linear operator, or alternatively, as
\begin{equation}\label{higher2}
    \omega = e^\eta d u + e^\eta h(d u),
\end{equation}
where $h$ is also a smooth linear operator, the coefficients of
which depend on $\Gamma;$ the latter variant yields better
regularity. In fact, the form $h(d u)$ is as smooth as $u,$ provided
$\Gamma$ is smooth.

A particular case of the above occurs when the coefficient $\Gamma$
is exact. In that case, the form $\omega$ satisfying (\ref{HF2}) is
said to be \emph{gradient recursive}. For gradient recursive
$k$-forms, (\ref{recursive}) assumes the simpler form (\ref{solu}),
in which $g$ is the identity.

\section{The linear case}

Corresponding to the physical example of Sec.\ 1.2, which
illustrates the nonlinear form of eqs.\ (\ref{HF1}), (\ref{HF2}), we
can illustrate the linear case by an even simpler physical example.

Condition (\ref{inte}) arises when a rigid body rotates in the
$xy$-plane at constant angular velocity $\tilde\omega.$ Taking the
axis of rotation to lie at the origin of coordinates, we write the
tangential velocity vector in the form
\[
\mathbf{v} = v_1\hat\imath + v_2\hat\jmath,
\]
where
\begin{equation}\label{rigid1}
    v_1 = -\tilde\omega y,\,\, v_2 = \tilde\omega x.
\end{equation}
Then
\begin{equation}\label{rigid2}
    \nabla\times\mathbf{v} = 2\tilde\omega\hat{k},
\end{equation}
so (\ref{inte}) is satisfied; see, \emph{e.g.}, Exercise 4.4 of
\cite{MTW}. In the sequel we take $\tilde\omega\equiv 1$ for
simplicity.

Equations (\ref{rigid1}) imply that
\[
\nabla\cdot\mathbf{v}=0,
\]
so we express the 1-form $\omega$ dual to $\mathbf{v}$ as a solution
of the \emph{linear} Hodge-Frobenius equations
\begin{eqnarray}
    \delta\omega=0,\nonumber\\
    d\omega=\Gamma\wedge\omega.\label{HF_lin}
\end{eqnarray}
Applying (\ref{solu}) and (\ref{grad_recurs}), we choose $\eta$ to
depend only on the distance $r$ from the axis of rotation. Then
\begin{eqnarray}
    d\eta\wedge\omega =\eta'(r)\cdot
r\,dxdy. \label{rigid2a}
\end{eqnarray}
In addition, eq.\ (\ref{rigid2}) implies that
\begin{equation}\label{rigid2b}
    d\eta\wedge\omega =d\omega = 2\,dxdy.
\end{equation}
Equating the right-hand sides of eqs.\ (\ref{rigid2a}) and
(\ref{rigid2b}), we conclude that $\eta(r)=2\log r$ and $\omega =
r^2du$ for $u(x,y) = \arctan\left(y/x\right).$ Then $|du|=r^{-1},$
so the singular structure of $u$ in the $xy$-plane is analogous to
that of the fundamental solution of Laplace's equation in
$\mathbb{R}^3.$ In particular, $u$ is singular at the origin of the
disc.

In this example, the Hodge-Frobenius equations themselves are only
defined on the punctured disc, as
\[
\Gamma\wedge\omega =
\left(\Gamma_1\omega_2-\Gamma_2\omega_1\right)dxdy = d\omega=2dxdy.
\]
Using (\ref{rigid1}), we can write this condition as an equation for
the inner product
\[
\left(%
\begin{array}{c}
  \Gamma_1 \\
  \Gamma_2 \\
\end{array}%
\right) \cdot \left(%
\begin{array}{c}
  x \\
  y \\
\end{array}%
\right) =2,
\]
which cannot be satisfied at the origin.

Thus singular solutions arise naturally in both the linear and
nonlinear Hodge-Frobenius equations.

Because eq.\ (\ref{recursive}) requires that the domain be
star-shaped, the conclusion that $\omega$ is representable as a
product $f\,du,$ where $f$ is nonvanishing, does not follow and is
in fact violated in our example, in which $f=r^2.$

We have presented a particularly simple model, in which the role of
condition (\ref{HF2}) is especially transparent. For more
sophisticated completely integrable models of rigid-body rotation
see, \emph{e.g.}, \cite{Bo} and references therein.

In the linear case we can accomplish easily what we cannot
accomplish at all in the nonlinear case: an integrability condition
sufficient to imply the smoothness of weak solutions.

\bigskip

\textbf{Proposition 1}. \emph{Let $\omega$ be a weak solution of the
linear Hodge-Frobenius equations (\ref{HF_lin}) on $\Omega.$ If
$|\Gamma|$ is bounded and $\omega\in L^p(\Omega)$ for $p>n,$ then
$\omega$ is continuous.}

\bigskip

\emph{Proof}. The Friedrichs mollification $\omega_h$ of $\omega$ is
a classical solution of (\ref{HF_lin}); see, \emph{e.g.}, Sec. 7.2
of \cite{GT}. Thus
\[
|d\omega_h|^p+|\delta\omega_h|^p+|\omega_h|^p\leq\left(|\Gamma|^p+1\right)|\omega_h|^p.
\]
Integrate and apply the $L^p$ Gaffney-G\"arding inequality (Lemma
4.7 of \cite{ISS}) to obtain
\[
||\nabla\omega_h||_p\leq C\left(\Gamma\right)||\omega_h||_p.
\]
Because $\omega\in L^p(\Omega),$ we can allow the mollification
parameter $h$ to tend to zero. The proof is completed by the Sobolev
Embedding Theorem.

\section{Geometric analogies}

\subsection{Hypersurface-orthogonal vector fields}

A unification of the two superficially different physical examples
of Secs.\ 1 and 3 can be found in their underlying geometry $-$ in
particular, in the relation of each vector field to hypersurfaces
created by the level sets of an associated scalar function $u.$ We
can think of $u$ informally as the potential for the conservative
field that would result from taking $\Gamma$ to be zero in eq.\
(\ref{HF2}). Whereas a conservative vector field is actually equal
to $\nabla u,$ the vector fields in Secs.\ 1 and 3 merely point in
the same direction as $\nabla u.$

A nonvanishing vector field $\mathbf{v}$ is said to be
\emph{hypersurface-orthogonal} whenever there exists a foliation of
hypersurfaces orthogonal to $\mathbf{v}.$ The foliated hypersurfaces
can be represented as level sets of a scalar function $u.$ That is,
one can write
\[
\mathbf{v} = \lambda(\mathbf{x})\nabla u
\]
for a nonvanishing function $\lambda.$ Conversely, a vector field
which can be written in this way is clearly hypersurface-orthogonal.
We conclude that a vector field $\mathbf{v}$ is
hypersurface-orthogonal if and only if the 1-form $\omega$ dual to
$\mathbf{v}$ satisfies
\[
\omega = \lambda du\;,
\]
with nonvanishing $\lambda$ (\emph{i.e.}, if and only if $\omega$
is completely integrable). Other equivalent conditions now follow
from the Frobenius Theorem.

Hypersurface-orthogonal vector fields arise naturally in general
relativity, particularly in connection with black-hole mechanics.
Introducing a tensor field
\[
\mathbf{B}_{\alpha\beta}\equiv\mathbf{v}_{\alpha;\beta},
\]
where the semi-colon denotes covariant differentiation with respect
to the spacetime metric connection, the condition that $\mathbf{v}$
be hypersurface-orthogonal implies that the antisymmetric part of
$\mathbf{B},$ called the \emph{rotation tensor}, vanishes (\cite{P},
Secs.\ 2.32, 2.33). For this reason, vector fields satisfying
(\ref{HF2}) are called \emph{rotation-free} in general relativity,
which is somewhat confusing in the context of the examples in Secs.
1 and 3.

\subsection{Twisted Born-Infeld equations}

A condition broadly analogous to (\ref{HF2}) arises if $\omega$ is a
Lie-algebra-valued 2-form satisfying the second Bianchi identity. In
that case, replacing $\omega$ by $F_A,$ where $A$ is a
Lie-algebra-valued 1-form, we have
\begin{equation}\label{bianchi}
    dF_A = -[A,F_A],
\end{equation}
where $[\,,\,]$ denotes the Lie bracket, an equation which resembles
(\ref{HF2}).

Precisely, let $X$ be a vector bundle over a smooth, finite,
oriented, \textit{n}-dimensional Riemannian manifold $M.$ Suppose
that $X$ has compact structure group $G\subset SO(m).$ Let $A\in
\Gamma \left( M,ad\,X\otimes T^{\ast }M\right) $ be a connection
1-form on $X$ having curvature 2-form
\[
F_{A}=dA+\frac{1}{2}\left[ A,\,A\right] =dA+A\wedge A,
\]
where [\ ,\ ] is the bracket of the Lie algebra $\Im ,$ the fiber of
the adjoint bundle $ad\,X.$ Sections of the automorphism bundle
$Aut\,X,$ called
\textit{gauge transformations,} act tensorially on $F_{A}$ but affinely on $%
A $; see, \textit{e.g.,} \cite{MM}.

Consider energy functionals having the form (\ref{energy}), where $Q=|F_{A}|^{2}=\left%
\langle F_{A},F_{A}\right\rangle $ is an inner product on the fibers
of the bundle $ad\,X\otimes \Lambda ^{2}\left( T\;^{\ast }M\right)
.$ The inner
product on $ad\,X$ is induced by the normalized trace inner product on $%
SO(m) $ and that on $\Lambda ^{2}\left( T\;^{\ast }M\right) ,$ by
the exterior product $\ast \left( F_{A}\wedge \ast F_{A}\right).$

A nonabelian variational problem analogous to eqs.\ (\ref{HF2}),
(\ref{HB}) is described briefly in Sec.\ 5.1 of \cite{O3}. One is
led to consider smooth variations taken in the infinitesimal
deformation space of the connection and having the explicit form
\begin{eqnarray}
    var\left( E\right) =\int_M\rho ( Q) var(Q) dM=
    \int_M\rho(Q)\frac{d}{dt}_{|t=0}|F_{A+t\psi}|^2dM \nonumber\\
    =\int_M\rho (Q) \frac{d}{dt}_{|t=0}\left|
F_{A}+tD_A\psi+t^{2}\psi\wedge \psi\right| ^{2}dM, \label{curvar}
\end{eqnarray}
where $D_A=d+[A,\,]$ is the exterior covariant derivative in the
bundle. The Euler-Lagrange equations are
\begin{equation}\label{nonab}
    \delta \left( \rho (Q)F_{A}\right) =-\ast \left[ A,\ast \rho
(Q)F_{A}\right] .
\end{equation}
In addition, we have the Bianchi identity (\ref{bianchi}).

Writing eq.\ (\ref{HF2}) in components
\[
d\omega^a=\Gamma^a_b\wedge \omega^a,
\]
we observe that if $-\Gamma$ is interpreted as a connection 1-form,
then (\ref{HF2}) can be interpreted as the vanishing of an exterior
covariant derivative, which is the content of eq.\ (\ref{bianchi}).
Moreover, the well known algebraic requirement that $\Gamma$ must
satisfy
\[
\left(d\Gamma^a_b-\Gamma^a_c\wedge\Gamma^c_b\right)\wedge \omega^b =
0
\]
(\emph{c.f.} eq.\ (4-2.3) of \cite{E}) is the zero-curvature
condition $\left[F_\Gamma,\omega\right]=0.$

If $G$ is abelian, then the Lie bracket vanishes and eqs.\
(\ref{nonab}) reduce to the system
\[
\delta \left\{ \rho \left[ Q(F_{A})\right] F_{A}\right\} =\delta
\left\{ \rho \left[ Q(dA)\right] dA\right\} =0,
\]
a nonlinear Hodge equation analogous to taking $\Gamma=0$ in eqs.\
(\ref{HF1}), (\ref{HF2}).

Equations (\ref{bianchi}) reduce in the abelian case to the
equations for the equality of mixed partial derivatives,
\[
d^{2}A=0.
\]

If $\rho\equiv 1,$ then eqs.\ (\ref{nonab}) are the Yang-Mills
equations, describing quantum fields in the classical limit. These
resemble a version of eq.\ (\ref{HB}) for 2-forms with $\rho\equiv
1,$ with the Bianchi identity (\ref{bianchi}) playing the role of
eq.\ (\ref{HF2}).

Whereas the Yang-Mills equations do not have the nonlinear structure
of (\ref{HB}) for non-constant $\rho,$ those of the Born-Infeld
model for electromagnetism are equivalent to (\ref{HB}) for
differential forms of degree 2 with
\begin{equation}\label{BIrho}
     \rho(Q)=\left(1+|F_A|^2\right)^{-1/2}.
\end{equation}
This model was introduced in \cite{BI} in order to produce a model
of electromagnetism that does not diverge when the source is a point
charge. Geometric aspects of the model are investigated in \cite{G}
and its analytic aspects, in \cite{Y}. A mathematical generalization
of the Born-Infeld model to nonabelian variational equations was
proposed in \cite{O1} and further studied in \cite{SSY}; see also
\cite{I} for a related problem. The equations assume the form
(\ref{bianchi}), (\ref{nonab}) for an appropriate choice of
$\rho(Q).$

Geometrically, the nonabelian model puts a twist in the principal
bundle corresponding to the configuration space of solutions. Thus
eqs.\ (\ref{bianchi}), (\ref{nonab}) are called \emph{twisted
nonlinear Hodge equations} \cite{O3}. A different approach to
generalizing nonlinear Hodge theory to bundle-valued connections,
which is based on the formulation of a natural class of
boundary-value problems, is introduced in \cite{Ma}; but interior
estimates would be required in order to extend the theory of
\cite{Ma} to twisted forms of the equations considered in this
paper. The derivation of such estimates is a goal of Sec.\ 5.2.

\section{Analysis}

We do not expect rotational fields of any kind to be very smooth. In
particular, assumption (\ref{HF2}) may produce caustics; see, for
example, the discussion in Sec.\ 2.4 of \cite{P}. However, it is
reasonable to seek conditions under which the field remains bounded
at a singularity or under which the field equations remain uniformly
elliptic.

Because these conditions will be derived for a large class of mass
densities, the strength of the estimates obtained will depend on the
integrability with respect to $\mathbf{x}$ of a given choice of mass
density $\rho\left(Q(\mathbf{x})\right).$

In the sequel we denote by $C$ generic positive constants, the value
of which may change from line to line. We follow an analogous
convention for continually updated \emph{small} positive constants
$\varepsilon.$ Repeated indices are to be summed from 1 to $n.$

\subsection{Variational structure}

The energy $E$ of the field $\omega$ on $\Omega$ is defined by eq.\
(\ref{energy}), where $Q = Q(\omega)$ is defined as in (\ref{speed})
for $\omega\in\Lambda^1$ given by (\ref{solu}). Then $\forall\psi\in
C_0^\infty(\Omega),$ the variations of $E$ are computed as
\[
var(E) = \frac{d}{dt}E\left(u+t\psi\right)_{|t=0} =
\frac{1}{2}\int_\Omega\rho(Q)\frac{d}{dt}Q\left(u+t\psi\right)_{|t=0}d\Omega
\]
\[
=\frac{1}{2}\int_\Omega\rho(Q)\frac{d}{dt}\left[e^{2\eta}|d\left(u+t\psi\right)|^2\right]_{|t=0}d\Omega
=
\]
\[
\frac{1}{2}\int_\Omega\rho(Q)e^{2\eta}\frac{d}{dt}\left(|du|^2+2\langle
du,td\psi\rangle+t^2|d\psi|^2\right)_{|t=0}d\Omega=
\]
\[
\int_\Omega\rho(Q)e^{2\eta}\langle du,d\psi\rangle
d\Omega=\int_\Omega\langle\rho(Q)e^{2\eta}du,d\psi\rangle d\Omega=
\]
\[
\int_\Omega d\langle\rho(Q)e^{2\eta}du,\psi\rangle
d\Omega+\int_\Omega\langle\delta\left[\rho(Q)e^\eta\omega\right],\psi\rangle
d\Omega=\int_\Omega\langle\delta\left[\rho(Q)e^\eta\omega\right],\psi\rangle
d\Omega,
\]
as $\psi$ has compact support in $\Omega.$ At a critical point,
$var(E)=0,$ or
\begin{equation}\label{aster}
    \delta\left[\rho(Q)e^\eta\omega\right]=0,
\end{equation}
so the variational formulation yields a ``weighted" form of the
continuity equation (\ref{HF1}). The presence of this weight adds an
inhomogeneous term to the un-weighted variant. To see this, we write
the local equation
\[
-\partial_i\left[\rho(Q)e^\eta\omega_i\right]=-\left(\partial_i\eta\right)e^\eta\rho(Q)\omega_i-e^\eta\partial_i\left[\rho(Q)\omega_i\right]=0.
\]
This local form corresponds to the invariant representation
\[
e^\eta\delta\left[\rho(Q)\omega\right]=e^\eta\rho(Q)\langle
d\eta,\omega\rangle.
\]
Thus we obtain eq.\ (\ref{HB}) for one-forms $\omega$
(\emph{i.e.}, for $k=1$) or equivalently, since
    $\delta\alpha=-\ast d\ast\alpha$, $\forall\alpha\in\Lambda^1$,
the equation
\begin{equation}
    d\ast\left[\rho(Q)\omega\right]=-d\eta\wedge\ast\left[\rho(Q)\omega\right].
\end{equation}

\subsubsection{Variational equations for $k$-forms}

Recall that an exterior differential form $\omega$ of degree $k$
which is recursive with  coefficient $\Gamma$ can be written in the
form (\ref{higher2}), where $u$ is a $\left(k-1\right)$-form (which
depends on $\omega$), and where the function $\eta$ and the linear
operator $h$ depend only on $\Gamma$ (Sec.\ 2.2).

We now compute the variation of the energy functional (\ref{energy})
among all forms $\omega + t\alpha$ satisfying
\[
d\left(\omega + t\alpha\right) = \Gamma \wedge \left(\omega +
t\alpha\right).
\]
Such forms satisfy (\ref{higher2}) with fixed $\eta$ and $h.$
Therefore the variations of $E$ are computed as
\[
var(E) = \frac{d}{dt}E\left(\omega +t\alpha\right)_{|t=0} =
\frac{1}{2}\int_\Omega\rho(Q)\frac{d}{dt}Q\left(\omega+t\alpha\right)_{|t=0}d\Omega
\]
\[
=\int_\Omega \langle \rho(Q) \omega, \alpha\rangle d\Omega,
\]
where $\alpha =  e^\eta \left(d v + h(dv)\right)$ for a
$\left(k-1\right)$-form $v$ which depends on $\alpha.$

Thus
\[
var(E) = \frac{1}{2}\int_\Omega \langle e^{\eta}\rho(Q)\omega , dv +
h(dv)\rangle d\Omega = \]
\[=\int_\Omega \langle e^{\eta}\rho(Q)\omega,
G dv\rangle d\Omega = \int_\Omega \langle G^T e^{\eta}\rho(Q)\omega,
dv\rangle d\Omega,
\]
where $G\equiv g_{ij} (\Gamma)$ is an $n\times n$ matrix and $G^T$
is its transpose.

As the forms $\alpha$ (and thus $v$) are assumed to have compact
support in $\Omega,$ setting $var(E)=0$ is equivalent to imposing
the condition
\[
\delta\left[G^T e^{\eta}\rho(Q)\omega\right]=0.
\]
Notice that if $\omega$ is gradient recursive, then $G$ is the
identity matrix and we recover (\ref{aster}); or equivalently, since
$\forall\alpha\in\Lambda^k$
\begin{equation}\label{formalad}
    \delta\alpha=(-1)^{nk+n+1}\ast d\ast\alpha,
\end{equation}
the equation
\begin{equation}\label{altHB}
    d\ast\left[\rho(Q)\omega\right]=-d\eta\wedge\ast\left[\rho(Q)\omega\right],
\end{equation}
\emph{i.e.}, eq.\ (\ref{HB}) for a gradient-recursive $k$-form
$\omega$ with coefficient $d\eta.$

The variations employed in this section, applied directly to a
2-form, are necessarily different from the variations of
(\ref{curvar}), which are applied instead to the Lie-algebra valued
connection 1-form $A.$

In the remainder of the paper we focus mainly on 1-forms. However,
many of the results will extend easily to gradient-recursive
$k$-forms. Moreover, the structure of expressions (\ref{higher1})
and (\ref{higher2}) suggest that, under appropriate technical
hypotheses on the linear operators $g$ or $h$, many of the results
will extend also to recursive $k$-forms with coefficient $\Gamma$
(not necessarily exact).

\subsection{When are solutions bounded at a singularity?}

Following \cite{U}, we find it convenient to introduce a function
$H(Q)$ which is defined so that
\begin{equation}\label{Hdef}
    H'(Q)=\frac{1}{2}\rho(Q)+Q\rho'(Q).
\end{equation}
Then ellipticity is equivalent to the condition that $H$ has
positive derivative with respect to $Q.$

In Theorem 7 and Corollary 8 of \cite{O1}, and in Theorem 6 and
Corollary 7 of \cite{O2}, $L^p$ conditions are derived which imply
the boundedness of solutions to eqs.\ (\ref{HF1}), (\ref{HF2}) on
domains that include singular sets of given co-dimension. We call
such theorems \emph{partial removable singularities theorems}, as
they imply that, although solutions may have jump discontinuities
at the singularity, they cannot blow up there. Those results
require the mass density $\rho$ to satisfy the inequality
\begin{equation}\label{nonnewt}
    C\left(K+Q\right)^q\leq
H'(Q)\leq C^{-1}\left(Q+K\right)^q
\end{equation}
for constants $q>0$ and $K\geq 0.$ (This hypothesis is imposed in
\cite{O2}, following Sec.\ 1 of \cite{U}; a somewhat stronger
hypothesis is imposed in \cite{O1}.)

There is obvious interest in deriving estimates for densities which
may not satisfy (\ref{nonnewt}). The focus of this section is to
obtain and exploit such estimates for a broad class of densities,
using the variational form of the equations. However, we retain the
condition that $\rho$ is positive, which is natural for
applications.

We impose an additional condition that arises from technical
considerations. If $\rho'(Q) > 0,$ we require that
\begin{equation}\label{hypo_pos}
    H'(Q)\leq C\rho(Q).
\end{equation}
(This inequality is satisfied automatically if $\rho'(Q)\leq 0.$) If
$\rho'(Q)<0$ we require instead that
\begin{equation}\label{hypo_neg}
    H'(Q)\geq C\rho(Q).
\end{equation}
(This inequality is satisfied automatically if $\rho'(Q)\geq 0.$)
Note that (\ref{hypo_neg}) implies  (\ref{subsonic}) under our
assumption $\rho(Q)>0.$ In the case $\rho'(Q)<0,$ eq.\ (\ref{HB}) is
uniformly elliptic (i.e. (\ref{unison}) is satisfied) whenever
(\ref{hypo_neg}) is satisfied and $\rho$ is \emph{noncavitating}:
bounded below away from zero. In the case $\rho'(Q)>0,$ condition
(\ref{unison}) is satisfied whenever $\rho(Q)$ is bounded above.

Densities which satisfy (\ref{nonnewt}) satisfy the hypotheses
(\ref{hypo_pos}), (\ref{hypo_neg}). However, there are many
densities which satisfy (\ref{hypo_pos}), (\ref{hypo_neg}) but do
not satisfy (\ref{nonnewt}). Among the latter are densities for
which the value of the exponent $q$ on the left-hand side of
inequality (\ref{nonnewt}) differs from its value on the right-hand
side, and certain densities for which the value of $q$ in
(\ref{nonnewt}) is negative.

As a simple illustration, consider the class of densities
\begin{equation}\label{nonnewt1}
    \rho(Q) =\left(K+Q\right)^q,\,-1/2<q<0,\, K>0.
\end{equation}
Such densities do not satisfy condition (\ref{nonnewt}) and cavitate
as $Q$ tends to infinity. They arise, for example, in connection
with models of pseudo-plastic non-Newtonian fluids \cite{B}. In this
section we will obtain sufficient conditions for a bound on
$Q\rho(Q)$ which is valid as $Q$ tends to infinity. Under additional
conditions, this bound extends to possibly singular solutions of the
system (\ref{HF2}), (\ref{HB}) with density given by
(\ref{nonnewt1}). Note that for such densities, a bound on the
product $Q \rho(Q)$ implies an asymptotic bound on the norm $Q$ of
the solution itself.

Because conditions (\ref{hypo_pos}) and (\ref{hypo_neg}) are used in
a crucial way to establish both the subelliptic estimates and the
ellipticity of the second-order operator, they appear to provide a
mathematically natural generalization of condition (\ref{nonnewt}).

\bigskip

\textbf{Lemma 2}. \emph{Let condition (\ref{hypo_neg}) be satisfied.
Then $H$ can be chosen so that}
\begin{equation}\label{hbound}
    Q\rho(Q) \leq CH(Q).
\end{equation}

\bigskip

\emph{Proof}. If $\rho'(Q)$ is non-negative, then (\ref{hbound}) is
always satisfied with $C=2$ if we choose $H(0)\geq 0.$ (Note that
assumption (\ref{hypo_neg}) is not needed in this case.) In order to
see this, let
\[
\Phi(Q) \equiv 2H(Q)-Q\rho(Q).
\]
Then $\Phi(0) = H(0)\geq 0,$ and
\[
\Phi'(Q) = 2H'(Q)-\rho(Q)-Q\rho'(Q)=Q\rho'(Q)\geq 0.
\]
Thus $\Phi(Q)$ remains nonnegative on the entire range of $Q.$

If $\rho'(Q)\leq 0,$ we assume (\ref{hypo_neg}). Then in particular,
\begin{equation}\label{hypo_neg_e}
    2H'(Q)\geq \varepsilon\rho(Q),
\end{equation}
where we take $\varepsilon$ to be so small that it lies in the
interval $\left(0,1\right).$ Inequality (\ref{hypo_neg_e}) can be
written in the form
\begin{equation}\label{vareps}
    \left(1-\varepsilon\right)\rho(Q)+2Q\rho'(Q)\geq 0.
\end{equation}
Define a constant $c$ by the formula
\[
c = \frac{1+\varepsilon}{2\varepsilon}.
\]
In terms of $c,$ (\ref{vareps}) can be written in the form
\[
\left(c-1\right)\rho(Q)+2\left(c-1/2\right)Q\rho'(Q)\geq 0.
\]
We can convert this expression into the differential inequality
\[
2cH'(s)\geq \frac{d}{ds}\left[s\rho(s)\right], \, s\in[0,Q].
\]
Integrate the inequality over $s,$ using $H(0)=0.$ We obtain
(\ref{hbound}).

\bigskip

\textbf{Lemma 3}. \emph{Let the 1-forms $\Gamma$ and $\omega$
smoothly satisfy (\ref{HF2}) and (\ref{HB}). Let $\rho>0$ satisfy
conditions (\ref{hypo_pos}) and (\ref{hypo_neg}). Then}
\begin{equation}\label{subell}
    \Delta H +(-1)^{3n}\nabla\cdot\left\{\ast
\left[\omega\wedge\ast\left(\rho'(Q)dQ\wedge\omega\right)\right]\right\}
+C\left(|\Gamma|^2+|\nabla\Gamma|\right) H \geq 0.
\end{equation}

\bigskip

\emph{Proof}. We have (\emph{c.f.} Sec.\ 1 of \cite{U})
\begin{eqnarray}
    \langle\omega,\Delta\left[\rho(Q)\omega\right]\rangle =
\partial_i\langle\omega,\partial_i\left(\rho(Q)\omega\right)\rangle-\langle\partial_i\omega,\partial_i\left(\rho(Q)\omega\right)\rangle\nonumber\\
=\Delta H(Q)- \left[\rho(Q) \langle
\partial_i\omega,\partial_i\omega\rangle+\rho'(Q)\langle\partial_i\omega,\omega\rangle\partial_i
Q\right],\label{kar}
\end{eqnarray}
where
\[
\Delta H(Q) =
\partial_i\left(\partial_i H(Q)\right)=\partial_i\left(H'(Q)\partial_i
Q\right).
\]
Writing $\partial_i Q = 2\langle\partial_i\omega,\omega\rangle,$ we
rewrite (\ref{kar}) in the form
\begin{equation}\label{SE1A}
    \langle\omega,\Delta\left[\rho(Q)\omega\right]\rangle =
\Delta
H(Q)-\rho(Q)|\nabla\omega|^2-2Q\rho'(Q)\left|d|\omega|\right|^2.
\end{equation}
Applying eq.\ (\ref{SE1A}) to the operator identity $\Delta =
-\left(\delta d + d\delta\right)$ and using (\ref{HB}), we write
\begin{eqnarray}\label{esti}
0=\langle\omega,\Delta\left[\rho(Q)\omega\right]\rangle
+\langle\omega,\delta
d\left(\rho(Q)\omega\right)\rangle+\langle\omega,d\delta\left(\rho(Q)\omega\right)\rangle\nonumber\\
=\Delta
H(Q)-\gamma+\langle\omega,\delta\left(d\rho\wedge\omega\right)\rangle+\tau_1+\tau_2,\label{below},
\end{eqnarray}
where
\[
\gamma=\rho(Q)|\nabla\omega|^2+2Q\rho'(Q)\left|d|\omega|\right|^2,
\]
\[
 \tau_1=\langle\omega,\delta
\left(\rho(Q)d\omega\right)\rangle,
\]
and
\[
\tau_2=\langle\omega,d\left[\rho(Q)\langle\Gamma,\omega\rangle\right]\rangle.
\]
Define
\[
L_\omega(H)\equiv\Delta H
+\langle\omega,\delta\left(d\rho(Q)\wedge\omega\right)\rangle.
\]
Then (\ref{esti}) can be written in the compact form
\begin{equation}\label{esti_comp}
    L_\omega H + \tau_1+\tau_2=\gamma.
\end{equation}
We have
\[
\tau_1=\langle\omega,\rho(Q)\delta\left(\Gamma\wedge\omega\right)-\langle
d\rho(Q),d\omega\rangle\rangle
\]
\[
\leq
|\omega|\left[\rho(Q)|\delta\left(\Gamma\wedge\omega\right)|+|d\rho(Q)||d\omega|\right]\leq
\]
\[
|\omega|\left[\rho(Q)|\nabla\Gamma||\omega|+\rho(Q)|\Gamma||\nabla\omega|+|\rho'(Q)dQ||d\omega|\right].
\]
Applying (\ref{HF2}) to the last term on the right-hand side and
using
\[
|dQ|=2|\omega||d|\omega||,
\]
we obtain
\begin{eqnarray}
\tau_1\leq
Q\rho(Q)|\nabla\Gamma|+\rho(Q)|\Gamma||\omega||\nabla\omega|+2Q|\rho'(Q)||\Gamma||\omega||d|\omega||\leq\nonumber\\
Q\rho(Q)|\nabla\Gamma|+\rho(Q)\left[\frac{Q|\Gamma|^2}{2\varepsilon}+\frac{\varepsilon}{2}|\nabla\omega|^2\right]\nonumber\\
+2Q|\rho'(Q)|\left(\frac{\varepsilon}{2}|d|\omega||^2+\frac{|\Gamma|^2Q}{2\varepsilon}\right).\label{tau1est}
\end{eqnarray}
Estimating $\tau_2$ yields the same terms:
\begin{eqnarray}
\tau_2\leq
Q\rho(Q)|\nabla\Gamma|+|\omega||\Gamma|\left[\rho(Q)|\nabla\omega|+|\omega|\nabla\rho(Q)|\right]=\nonumber\\
Q\rho(Q)|\nabla\Gamma|+|\omega||\Gamma|\rho(Q)|\nabla\omega|+2Q|\Gamma||\rho'(Q)||\omega||d|\omega||,\label{tau2}
\end{eqnarray}
which is bounded by the right-hand side of (\ref{tau1est}). That is,
\begin{eqnarray}
\tau_1+\tau_2\leq2Q\rho(Q)|\nabla\Gamma|+
\frac{Q}{\varepsilon}\left[\rho(Q)+2Q|\rho'(Q)|\right]|\Gamma|^2+\nonumber\\
\varepsilon\rho(Q)|\nabla\omega|^2+2Q\varepsilon|\rho'(Q)||d|\omega||^2.\label{tau12}
\end{eqnarray}
Applying inequality (\ref{tau12}) to eq.\ (\ref{esti_comp}), we
obtain, in the case $\rho'(Q)>0,$ the estimate
\begin{eqnarray}
L_\omega(H)+2\left[Q\rho(Q)|\nabla\Gamma|+\frac{QH'(Q)}{\varepsilon}|\Gamma|^2\right]\geq
\nonumber\\
\left(1-\varepsilon\right)\left[\rho(Q)|\nabla\omega|^2+2Q\rho'(Q)|d|\omega||^2\right]\geq
\left(1-\varepsilon\right)H'(Q)|d|\omega||^2,\label{almost}
\end{eqnarray}
the inequality on the right following from Kato's inequality. We
apply Lemma 2 and (\ref{hypo_pos}) to terms on the extreme left-hand
side of inequality (\ref{almost}):
\[
2\left[Q\rho(Q)|\nabla\Gamma|+\frac{QH'(Q)}{\varepsilon}|\Gamma|^2\right]\leq
\]
\[
C\left[H(Q)|\nabla\Gamma|+Q\rho(Q)|\Gamma|^2\right]\leq
\]
\[
C\left(|\nabla\Gamma|+|\Gamma|^2\right)H(Q).
\]
We now have, for the case $\rho'(Q)>0,$ the estimate
\begin{equation}\label{ellest}
L_\omega(H)+C\left(|\nabla\Gamma|+|\Gamma|^2\right)H\geq
\left(1-\varepsilon\right)H'(Q)|d|\omega|^2\geq 0.
\end{equation}
In the case $\rho'(Q)<0,$ we also apply (\ref{tau12}) to
(\ref{esti_comp}); but here we obtain
\begin{eqnarray}
L_\omega(H)+2\left[Q\rho(Q)|\nabla\Gamma|+\frac{QH'(Q)}{\varepsilon}|\Gamma|^2\right]\geq\nonumber\\
\left(1-\varepsilon\right)\left[\rho(Q)|\nabla\omega|^2+2Q\rho'(Q)|d|\omega||^2\right]+4\varepsilon
Q\rho'(Q)|d|\omega||^2.\label{neg_case}
\end{eqnarray}
In the case $\rho'(Q)<0,$ condition (\ref{hypo_neg}) yields
\[
-Q\rho'(Q)< \frac{1}{2}\rho(Q)\leq CH'(Q),
\]
thus
\begin{equation}\label{lemma}
    Q|\rho'(Q)|\leq CH'(Q).
\end{equation}
(Note that (\ref{lemma}) is automatic for $\rho'(Q)\geq 0.$)
Applying Kato's inequality and (\ref{lemma}) to (\ref{neg_case})
yields
\[
L_\omega(H)+2\left[Q\rho(Q)|\nabla\Gamma|+\frac{QH'(Q)}{\varepsilon}|\Gamma|^2\right]\geq
\]
\[
\left(1-\varepsilon-4\varepsilon C\right)H'(Q)|d|\omega||^2\geq 0.
\]
Applying (\ref{hypo_pos}) and (\ref{hbound}) to the left-hand side
of this inequality, we obtain (\ref{ellest}) (for updated
$\varepsilon$) for the case $\rho'(Q)<0$ as well.

It remains only show that the operator $L_\omega(H)$ can be put into
divergence form, at the cost of absorbing another lower-order term.
The second term of $ L_\omega(H)$ can be written in the form
\[
\langle\omega,\delta\left(d\rho(Q)\wedge\omega\right)\rangle=\ast\left[\omega\wedge\ast\delta
\left(d\rho(Q)\wedge\omega\right)\right] =
\]
\[
\ast
d\left[\omega\wedge\ast\left(d\rho(Q)\wedge\omega\right)\right]-\ast\left[d\omega\wedge\ast
\left(d\rho(Q)\wedge\omega\right)\right],
\]
where
\[
-\ast\left[d\omega\wedge\ast
\left(d\rho(Q)\wedge\omega\right)\right]=-\ast\left[\Gamma\wedge\omega\wedge\ast
\left(d\rho(Q)\wedge\omega\right)\right]\geq
\]
\[
-2Q|\Gamma||\rho'(Q)||\omega||d|\omega||,
\]
which is estimated in the same way as the last term in the sum on
the extreme right-hand side of (\ref{tau2}).

Taking into account (\ref{formalad}) we obtain, for any $k$-form
$\alpha,$
\begin{eqnarray}
    \ast d\alpha = (-1)^{k\left(n-k\right)}\ast
d\left(\ast\ast\right)\alpha=\nonumber\\
\left(-1\right)^{k\left(n-k\right)}\left(\ast
d\ast\right)\ast\alpha=
\left(-1\right)^{2kn+n+1-k^2}\delta\ast\alpha.\label{diverg}
\end{eqnarray}
Taking $k=1$ and
\[
\alpha =\ast
\left[\omega\wedge\ast\left(\rho'(Q)dQ\wedge\omega\right)\right]
\]
in eq.\ (\ref{diverg}), we can express the operator $\delta$ in that
equation as a divergence. This allows us to write (\ref{ellest})
(again updating $C$ and $\varepsilon$) in the form
\[
\Delta H +(-1)^{3n}\nabla\cdot\left\{\ast
\left[\omega\wedge\ast\left(\rho'(Q)dQ\wedge\omega\right)\right]\right\}
+C\left(|\Gamma|^2+|\nabla\Gamma|\right) H \geq
\]
\[
\left(1-\varepsilon\right)H'(Q)|d|\omega||^2\geq 0.
\]
This completes the proof of Lemma 3.

\bigskip

\textbf{Lemma 4}. \emph{Under the hypotheses of Lemma 3, the
operator}
\[
\mathcal L_\omega(H) \equiv \Delta H
+(-1)^{3n}\nabla\cdot\left\{\ast
\left[\omega\wedge\ast\left(\rho'(Q)dQ\wedge\omega\right)\right]\right\}
\]
\emph{is a uniformly elliptic operator on} $H.$

\bigskip

\emph{Proof}. We argue as in Sec.\ 1 of \cite{U} and Sec.\ 4 of
\cite{O2}, but without using hypothesis (\ref{nonnewt}).

Define a map $\beta_\omega:\Lambda^0\rightarrow\Lambda^{k+1}$ by the
explicit formula
\begin{equation}\label{KU_alt}
    \beta_\omega:\mu\rightarrow d\mu\wedge\omega,
\end{equation}
for $\mu\in\Lambda^0$ and $\omega\in\Lambda^k.$ Then we can write
the variational form of the Hodge-Frobenius equations (\ref{HF2}),
(\ref{solu}), (\ref{HB}) in the alternate form
\begin{equation}\label{HF2A}
    d\sigma_i=-\beta_{\sigma_i}(\eta),\,\,\sigma=1,2,
\end{equation}
where $\sigma_1=\ast\rho(Q)\omega$ and $\sigma_2=-\omega.$ The
``irrotational" case $d\omega=0$ can be recovered as the special
case of (\ref{HF2A}) in which
\[
\beta_\omega(\eta) = d\left(\eta\omega\right),
\]
in which case eq.\ (\ref{KU_alt}) implies that
\[
d\eta\wedge\omega=d\eta\wedge\omega+\eta d\omega
\]
(\emph{c.f.} Sec.\ 4 of \cite{O2} and Sec.\ 1 of \cite{U}).

Moreover, writing
\[
\beta_\omega(g)=dg\wedge\omega
\]
for some 0-form $g,$ we compute for arbitrary compactly supported
$\mu\in\Lambda^{k+1},$
\[
\langle\mu, dg\wedge\omega\rangle =
\ast\left(dg\wedge\left(\ast\ast\right)\left(\omega\wedge\ast\mu\right)\right)
=\langle dg,\ast\left(\omega\wedge\ast\mu\right)\rangle = \langle
g,\delta\ast\left(\omega\wedge\ast\mu\right)\rangle.
\]
So the map $\beta_\omega^\ast:\Lambda^{k+1}\rightarrow\Lambda^0$
defined by the explicit formula
\[
\beta_\omega^\ast(\mu) = \delta\ast\left(\omega\wedge\ast\mu\right)
\]
is the formal adjoint of $\beta_\omega.$

In terms of the maps $\beta_\omega$ and $\beta_\omega^\ast,$ we can
write
\[
(-1)^{n+1}\nabla\cdot\left\{\ast
\left[\omega\wedge\ast\left(\rho'(Q)dQ\wedge\omega\right)\right]\right\}=\beta^\ast_\omega\beta_\omega[\rho]
\]
\[
=\beta^\ast_\omega\left[\mu_\omega(H)\right]
\]
for $\mu_\omega$ satisfying
\[
    \mu_\omega(H)=\frac{\rho'(Q)}{H'(Q)}dH\wedge\omega.
\]
Using $\beta_\omega^\ast,$ we write the inequality of Lemma 3 in the
form
\[
\mathcal L_\omega(H) + \mbox{ lower-order terms }\geq 0,
\]
where
\[
\mathcal L_\omega(H) = \Delta H -
\beta^\ast_\omega\left[\mu_\omega(H)\right].
\]
Writing
\[
\mathcal L_\omega(H) =
\partial_k\left(\alpha^{jk}\partial_j\right)H,
\]
we find that if $\rho'(Q)<0,$ then (\ref{lemma}) implies that
\[
1\leq \alpha_{kj}+\frac{Q|\rho'(Q)|}{H'(Q)}\leq \alpha_{kj}+C.
\]
If $\rho'(Q)\geq 0,$ then we write
\[
\nabla\cdot\left\{\left[1-\frac{Q\rho'(Q)}{H'(Q)}\right]\nabla
H\right\}=\nabla\cdot\left\{\left[\frac{\rho(Q)}{2H'(Q)}\right]\nabla
H\right\};
\]
condition (\ref{hypo_pos}) implies that there is a positive constant
$c$ such that
\[
\frac{c}{2}\leq\frac{\rho(Q)}{2H'(Q)}\leq 1;
\]
$c$ is the reciprocal of the constant $C$ in (\ref{hypo_pos}). This
completes the proof of Lemma 4.

\bigskip

These three lemmas easily yield:

\bigskip

\textbf{Theorem 5}. \emph{Under the hypotheses of Lemma 3, the
product $Q\rho(Q)$ is locally bounded above by the $L^2$-norm of
$H.$}

\bigskip

The proof of Theorem 5 is given in Sec.\ 7.1.

\bigskip

Of course the integrability of $H$ depends on $\rho.$ But for any
given $\rho$ in $C^1(Q),$ $H$ can be computed explicitly by
integrating (\ref{Hdef}). Even if $\rho$ cavitates, Theorem 5 yields
asymptotic information about the fastest rate at which $Q$ can blow
up. Nevertheless, Theorem 5 is ultimately not very useful, due to
the hypothesis that the solutions are smooth. It would become more
useful if it could be applied to singular solutions. We will find
that the partial removable singularities theorems proven in
\cite{O2} for singular sets of prescribed codimension extend to our
conditions on $\rho$ under slightly different hypotheses.

Initially, we treat the special case of an isolated point
singularity, for which the proof is somewhat simpler than the proof
for higher-order singularities and the range of applicable
dimensions somewhat larger.

\bigskip

\textbf{Theorem 6}. \emph{Let the hypotheses of Lemma 3 be satisfied
on $\Omega\backslash\{p\},$ where $p$ is a point of $\mathbb{R}^n$
and $n>2.$ If $H\in L^{2n/\left(n-2\right)}(\Omega)$ and if the
function}
\begin{equation}\label{eff}
    f\equiv|\nabla\Gamma|+|\Gamma|^2
\end{equation}
\emph{is sufficiently small in $L^{n/2}(\Omega),$ then $H$ is an
$H^{1,2}$-weak solution in a neighborhood of the singularity.}

\bigskip

The proof of Theorem 6 is given in Sec.\ 7.2.

\bigskip

\textbf{Corollary 7}. \emph{Let the hypotheses of Theorem 6 be
satisfied and, in addition, let the function $f$ given by
(\ref{eff}) satisfy the growth condition}
\begin{equation}\label{mor}
    \int_{B_r\left(\mathbf{x}_0\right)\cap\Omega}|f|^{n/2}d\Omega\leq
    Cr^\kappa
\end{equation}
\emph{for some $\kappa >0,$ where $B_r\left(\mathbf{x}_0\right)$ is
an $n$-disc of radius $r,$ centered at $\mathbf{x}_0.$ Then the
conclusion of Theorem 5 remains valid.}

\bigskip

\emph{Proof}. Apply Theorem 5.3.1 of \cite{M} to the conclusion of
Theorem 6, following the proof of Theorem 5.

\bigskip

Note that the singularity in Corollary 7 is in the solution,
rather than in the underlying metric (\emph{c.f.} \cite{Sm}, in
which metric point singularities are considered in the case
$\Gamma\equiv 0$). Corollary 7 extends to higher-order
singularities in spaces of sufficiently high dimension; but the
proof requires more delicate test functions.

\bigskip

\textbf{Theorem 8}. \emph{Let the pair $\omega,$ $\Gamma$ smoothly
satisfy eqs.\ (\ref{HF2}), (\ref{HB}), with $\rho$ satisfying the
hypotheses of Lemma 3, on $\Omega\backslash\Sigma.$ Here $\Sigma$
is a compact singular set, of dimension $0\leq m<n-4,$ completely
contained in a sufficiently small $n$-disc $D$ which is itself
completely contained in the interior of $\Omega.$ Let $H(Q)$ lie
in $L^{2\beta\gamma_1}(D)\cap L^{2\gamma_2}(D),$ where $\beta =
\left(n-m-\varepsilon\right)/\left(n-m-2-\varepsilon\right)$ for
$1/2<\gamma_1<\gamma_2.$ If the function $f$ given by (\ref{eff})
satisfies the growth condition (\ref{mor}), then the conclusion of
Theorem 5 remains valid.}

\bigskip

The proof of Theorem 8 is given in Sec.\ 7.3.

\bigskip

The conditions imposed on $\rho$ in this section also lead to
extensions of known results for the conventional case $\Gamma\equiv
0.$ In particular, we consider equations which, expressed in
components, have the weak form
\begin{equation}\label{weakder}
    \int_\Omega\left[\rho(Q)u_{x_k}\right]_{x_i}\varphi_{x_k}d\Omega =
0,
\end{equation}
where $\varphi\in C_0^\infty(\Omega)$ is arbitrary and $Q=|du|^2.$
Such equations have been intensively studied in cases for which
$\rho\left(|du|^2\right)$ grows as a power of $|du|;$ \emph{c.f.}
eq.\ (3.10) of \cite{D}. Equation (\ref{weakder}) can be interpreted
as a weak derivative with respect to $x_i$ of the system
(\ref{HF1}), (\ref{HF2}) with $\Gamma\equiv 0.$

Define
\[
\mathfrak{H}(Q) \equiv Q\rho^2(Q).
\]
The system (\ref{HF1}), (\ref{HF2}) is elliptic precisely when
$\mathfrak{H}'(Q)>0.$ The following result extends Theorem 1 of
\cite{O2}, which requires the derivative of $\rho$ with respect to
$Q$ to be nonpositive. If $\rho$ is specified to be nonincreasing in
$Q,$ then the integrability of $\mathfrak{H}(Q)$ can be shown to
follow from finite energy. In the general case, we impose this
integrability as an independent hypothesis on the weak solution.

\bigskip

\textbf{Theorem 9}. \emph{Let the scalar function
$u\left(\mathbf{x}\right)$ satisfy equation (\ref{weakder}) with
$\rho$ bounded, positive, and noncavitating, and with
$\mathfrak{H}(Q)$ integrable. Assume conditions (\ref{subsonic}) and
(\ref{hypo_pos}). Then for every $n$-disc $D_R$ of radius $R$
completely contained in $\Omega$ there is a positive number
$\delta>0$ such that}
\[
\sup_{Q\in D_{\left(1-\delta\right)R}}\mathfrak{H}(Q)\leq
CR^{-n}\int_{D_R}\mathfrak{H}(Q)\ast 1,
\]
\emph{where $C$ depends on $\rho$ and $\delta$ but not on $Q$ or
$R.$}

\bigskip

\emph{Proof}. Define $Q_\ast$ to be the (possibly infinite) value of
$Q$ for which
\[
\sup_{Q\in\Omega}\rho(Q)=\rho\left(Q_\ast\right).
\]
Then
\begin{equation}\label{mult_est}
    \frac{\mathfrak{H}'(Q)}{\rho^2\left(Q_\ast\right)}
    =\frac{2\rho(Q)H'(Q)}{\rho^2\left(Q_\ast\right)}\leq
    C\left(\frac{\rho(Q)}{\rho\left(Q_\ast\right)}\right)^2\leq
    C,
\end{equation}
where $C$ is the constant of (\ref{hypo_pos}).

As in \cite{O2}, we initially estimate smooth solutions, and
subsequently extend the result to weak solutions, recovering the
derivatives as limits of difference quotients.

We choose test functions $\varphi^i$ having the form
\[
\varphi^i(\mathbf{x}) = u_{x_i}\tilde\mathfrak{H}^{\alpha/2}\zeta^2,
\]
where $\alpha >0,$ $\zeta\left(\mathbf{x}\right) \in
C_0^\infty\left(D_R\right),$ and
\begin{equation}\label{stable}
    \tilde\mathfrak{H} \equiv \mathfrak{H}(Q) + \varepsilon
\end{equation}
for a small positive parameter $\varepsilon$ (\emph{c.f.} \cite{D},
Sec.\ 3). Then
\[
\left[\rho(Q)u_{x_k}\right]_{x_i}\varphi^i_{x_k}=
\]
\[
2H'(Q)\left(u_{x_ix_k}\right)^2\tilde\mathfrak{H}^{\alpha/2}\zeta^2+\alpha\rho(Q)\left[H'(Q)\right]^2|\nabla
Q|^2\tilde\mathfrak{H}^{\left(\alpha-2\right)/2}\zeta^2
\]
\[
+2H'(Q)Q_{x_k}\tilde\mathfrak {H}^{\alpha/2}\zeta\zeta_{x_k} \equiv
i_1+i_2+i_3.
\]
Here
\[
i_1 \geq
2H'(Q)u_{x_ix_k}u_{x_ix_k}Q\rho^2(Q)\tilde\mathfrak{H}^{\left(\alpha-2\right)/2}\zeta^2=
\]
\[
\frac{1}{2}\rho^2(Q)H'(Q)|\nabla
Q|^2\tilde\mathfrak{H}^{\left(\alpha-2\right)/2}\zeta^2 \geq
C|\nabla\left(\tilde\mathfrak{H}^{\left(\alpha+2\right)/4}\right)|^2\zeta^2,
\]
where the last inequality follows from (\ref{mult_est}) and the
constant $C$ depends on $\alpha,$ $\rho^{-1}\left(Q_\ast\right),$
the constant of (\ref{mult_est}), and the lower bound of $\rho(Q).$

Similarly,
\[
i_2\geq\frac{\rho(Q)}{\rho\left(Q_\ast\right)}i_2\geq C
|\nabla\left(\tilde\mathfrak{H}^{\left(\alpha+2\right)/4}\right)|^2\zeta^2
\]
and
\[
i_3 = 2H'(Q)
Q_{x_j}\tilde\mathfrak{H}^{\alpha/2}\zeta\left(\mathbf{x}\right)\zeta_{x_j}.
\]
The latter quantity can be estimated by applying, as in the proof of
Lemma 3, the elementary algebraic inequality (\emph{Young's
inequality})
\[
ab\geq
-\left(\frac{\tilde\varepsilon}{2}a^2+\frac{1}{2\tilde\varepsilon}b^2\right),
\]
now taking
\[
a=2H'(Q)\tilde\mathfrak{H}^{\left(\alpha-2\right)/4}\zeta
Q_{x_j},\,\, b =
\tilde\mathfrak{H}^{\left(\alpha+2\right)/4}\zeta_{x_j},
\]
and $\tilde\varepsilon = \tilde\delta\rho^2(Q)$ for small
$\tilde\delta>0.$

Putting these estimates together, we obtain
\[
\int_\Omega|\nabla\left(\tilde\mathfrak{H}^{\left(\alpha+2\right)/4}\right)|^2\zeta^2d\Omega\leq
C\int_\Omega\tilde\mathfrak{H}^{\left(\alpha+2\right)/2}|\nabla\zeta|^2d\Omega.
\]
The proof for the smooth case is completed by applying the Moser
iteration as in expressions (9.5.8)-(9.5.12) of \cite{LU} and
subsequently letting the parameter $\varepsilon$ in (\ref{stable})
tend to zero. The proof is extended to the general case by applying
the difference-quotient method as in eq.\ (13) of \cite{O2}; see
also Lemma 2 of \cite{U}.

\section{Hodge-B\"acklund transformations}

Different choices of $\rho$ may sometimes be related by a special
kind of B\"acklund transformation which is based on properties of
the Hodge involution. We call these transformations
\emph{Hodge-B\"acklund}. Although this term does not seem to have
been used up to now, such generalized B\"acklund transformations
have a long history in diverse fields of mathematical physics. Our
aim in this section is to unify these various transformations, place
them in an invariant context, and extend them to the completely
integrable case.

Historically, the term \emph{B\"acklund transformation} has been
defined in many ways; see, \emph{e.g.}, \cite{RS} for the classical
theory. In the sequel we will use it in the general sense of a
function that maps a solution $a$ of a differential equation $A$
into a solution $b$ of a differential equation $B$ and vice-versa,
where $B$ may equal $A$ but $b$ will not equal $a.$

\subsection{Transformation of the Chaplygin mass density}

The mass density for the adiabatic and isentropic subsonic flow of
an ideal fluid has the form
\begin{equation}\label{gasdyn}
    \rho(Q) = \left(1-\frac{\gamma-1}{2}Q\right)^{1/(\gamma-1)},
\end{equation}
for $Q\in\left[0,2/\left(\gamma+1\right)\right),$  where $\gamma$ is
the \emph{adiabatic constant}: the ratio of specific heats for the
gas. The adiabatic constant for air is 1.4. Choosing $\gamma$ to be
2 we obtain, by an independent physical argument originally
introduced for one space dimension in \cite{R}, the mass density for
shallow hydrodynamic flow in the tranquil regime (\emph{c.f.} eq.\
(10.12.5) of \cite{St}). If we choose $\gamma$ to be $-1$ (a
physically impossible choice), we obtain the density of the minimal
surface equation \cite{K}, \cite{SS2}
\begin{equation}\label{minsurf}
    \rho(Q)=\frac{1}{\sqrt{1+Q}}
\end{equation}
Flow governed by this density is called \emph{Chaplygin flow}.
Despite the fact that the numbers $-1$ and 1.4 are not particularly
close, this choice of mass density has many attractive properties as
an approximation for (\ref{gasdyn}); see, \emph{e.g.}, \cite{DC} and
Ch.\ 5 of \cite{Be2}. These properties are, in general, retained in
the case of completely integrable flow described in Sec.\ 1.2.

If $\Gamma\equiv 0,$ eqs.\ (\ref{HF1}) and (\ref{HF2}) with
$\rho(Q)$ given by (\ref{minsurf}) describe, for $k=1,$
nonparametric minimal surfaces embedded in Euclidean space. If $k=2$
they describe electromagnetic fields in the Born-Infeld model, as in
(\ref{BIrho}).

More generally, we have the following result, which extends an
argument introduced for the case $\Gamma\equiv 0$ by Yang in
\cite{Y}; see also \cite{AP} and Theorem 2.1 of \cite{SSY}.

\bigskip

\textbf{Theorem 10}. \emph{Let the 1-form $\omega$ satisfy eqs.
(\ref{HF2}) and (\ref{HB}), with $\rho$ satisfying (\ref{minsurf}).
Then there exists an $\left(n-1\right)$-form $\xi$ with $|\xi|<1,$
satisfying equations analogous to (\ref{HF2}) and (\ref{HB}), but
with $\Gamma\equiv d\eta$ replaced by $\hat\Gamma\equiv
d\hat\eta=-d\eta$ and $\rho(Q)$ replaced by
\begin{equation}\label{MSLH}
    \hat\rho(|\xi|^2)\equiv\frac{1}{\sqrt{1-|\xi|^2}}\;.
\end{equation}
}

\bigskip

\emph{Proof}. Equation (\ref{HB}) can be interpreted as the
assertion that the $\left(n-1\right)$-form
\begin{equation}
\label{xi}
    \xi = \ast\left[\rho(Q)\omega\right]= \ast\left[\frac{\omega}{\sqrt{1+|\omega|^2}}\right]
\end{equation}
satisfies
\begin{equation}\label{fivealt}
    d\xi = d\hat\eta\wedge\xi,
\end{equation}
that is, equation (\ref{HF2}) with $\eta$ replaced by
$\hat\eta\equiv -\eta$. As a consequence, we conclude that this
$\left(n-1\right)$-form is also gradient recursive and, on domains
with trivial de Rham cohomology, there exists an
$\left(n-2\right)$-form $\sigma$ such that $\xi =
e^{\hat\eta}d\sigma$ (\emph{c.f.} Sec.\ 2.2, final paragraph).

Because the Hodge involution is an isometry,
\[
|\xi|^2=\frac{|\omega|^2}{1+|\omega|^2}
\]
or, equivalently,
\begin{equation}\label{iso}
\label{isometry} 1-|\xi|^2=\frac{1}{1+|\omega|^2}\;.
\end{equation} Note that eq.\ (\ref{iso}) implies $|\xi|^2<1$, as well as
\[
\rho(|\omega|^2)\hat\rho(|\xi|^2)=\frac{1}{\sqrt{1+|\omega|^2}}
\frac{1}{\sqrt{1-|\xi|^2}}= 1.
\]
This, together with \ (\ref{xi}), yields directly
\[
\ast\hat\rho(|\xi|^2)\xi = \ast^2
\rho(|\omega|^2)\hat\rho(|\xi|^2)\omega = (-1)^{n-1}\omega.
\]
Hence,
\[
d\ast (\hat\rho(|\xi|^2)\xi) = (-1)^{n-1}d\omega =
(-1)^{n-1}d\eta\wedge\omega =d\eta\wedge \ast
\left(\hat\rho(|\xi|^2)\xi\right)= -d\hat\eta\wedge \ast
\left(\hat\rho(|\xi|^2)\xi\right),
\]
which is equivalent to equation (\ref{HB}) for the gradient
recursive $\left(n-1\right)$-form $\xi$ with coefficient
$d\hat\eta$, where $\rho$ has been replaced by $\hat\rho$ (see
also (\ref{altHB})) . This completes the proof of Theorem 10.

\bigskip

The above argument carries over to any pairing of functions $\rho
(|\omega|^2)$, $\hat\rho(|\xi|^2)$, as long as their product is 1.
Moreover, it carries over essentially unchanged to
gradient-recursive $k$-forms $\omega$, as these would automatically
yield gradient-recursive $\left(n-k\right)$-forms $\xi$.

The same argument, with small modifications, extends Theorem 10 to
general (non-gradient-recursive) $k$-forms, in which we write
equations (\ref{HF2}) and (\ref{HB}) in terms of $\Gamma$ and
$\hat\Gamma$ rather than in terms of $\eta$ and $\hat{\eta}.$

In addition to the original, ``irrotational" version of Theorem 10
introduced in \cite{Y}, other aspects of the duality of mass
densities for nonlinear Hodge equations are presented in
\cite{SS1}, and in Sec.\ 4.2 of \cite{O3}. (The irrotational case
of the above argument is reviewed in Sec. 4.1 of \cite{O3}; note
the recurring misprint in the two paragraphs following eq.\ (35)
of that reference: $d\omega$ should be $du.$) Densities of the
form (\ref{MSLH}) arise in the study of maximal space-like
hypersurfaces \cite{C} and, in a completely different way,
harmonic diffeomorphisms \cite{T}.

\subsection{Transformation of the complex eikonal equation (\emph{after Magnanini and Talenti})}

As an example of the diverse fields in which these transformations
arise, and of the simplifying and unifying role of the Hodge
involution, we describe an example from complex optics. The
description of the local B\"acklund transformations follows the
analysis of Magnanini and Talenti \cite{MT1}, in which these
transformations were introduced; see also \cite{MT2} and \cite{MT3}.
We then reproduce the argument of \cite{MT1} in a simpler, invariant
context using the Hodge operator.

The \emph{eikonal equation} in $\mathbb{R}^2$ can be written in the
form
\begin{equation}\label{eik}
    \psi_x^2+\psi_y^2+\nu^2=0,
\end{equation}
where $\nu\left(x,y\right)$ is a given real-valued function. If we
write the solution $\psi\left(x,y\right)$ as a complex function
having the form
\[
\psi\left(x,y\right) = u\left(x,y\right)+iv\left(x,y\right)
\]
for real-valued functions $u$ and $v,$ then (\ref{eik}) is
equivalent to the first-order system
\begin{equation}\label{diff}
    u_x^2+u_y^2-v_x^2-v_y^2+\nu^2 =0,
\end{equation}

\begin{equation}\label{orthog}
    u_xv_x+u_yv_y=0.
\end{equation}

The function $\nu$ corresponds physically to the refractive index of
the medium through which the wavefront represented by the function
$\psi$ propagates. If
\begin{equation}\label{poshyp1}
    u_x^2+u_y^2>0,
\end{equation}
then $v_x$ and $v_y$ cannot both vanish, by eq.\ (\ref{diff}). In
that case, eq.\ (\ref{orthog}) can be expressed, in the language of
proportions, by the assertion that either
\begin{equation}\label{propo}
    \left[%
\begin{array}{c}
  v_x \\
  v_y \\
\end{array}%
\right] : \sqrt{v_x^2+v_y^2}=\pm\left[%
\begin{array}{c}
  -u_y \\
  u_x \\
\end{array}%
\right] : \sqrt{u_x^2+u_y^2},
\end{equation}
or
\begin{equation}\label{vanish}
    u_x=u_y=0.
\end{equation}

Treating these relations as a coupled system of scalar equations, we
have, under the same hypothesis,
\[
\frac{v_x}{\sqrt{v_x^2+v_y^2}}=\pm\frac{-u_y}{\sqrt{u_x^2+u_y^2}}
\]
and, using eq.\ (\ref{diff}),
\begin{equation}\label{exact1}
    v_x=\pm \left(-u_y\right)\sqrt{\frac{\nu^2}{u_x^2+u_y^2}+1}.
\end{equation}

Similarly, the scalar equation
\[
\frac{v_y}{\sqrt{v_x^2+v_y^2}}=\pm\frac{u_x}{\sqrt{u_x^2+u_y^2}},
\]
which also follows from (\ref{propo}), implies by analogous
operations the equation
\begin{equation}\label{exact2}
    v_y=\pm u_x\sqrt{\frac{\nu^2}{u_x^2+u_y^2}+1}.
\end{equation}

We can write the coupled system (\ref{exact1}), (\ref{exact2}) as a
vector equation of the form
\begin{equation}\label{vec1}
    \left[%
\begin{array}{c}
  v_x \\
  v_y \\
\end{array}%
\right]=\pm\sqrt{1+\frac{\nu^2}{u_x^2+u_y^2}}\left[%
\begin{array}{c}
  -u_y \\
   u_x \\
\end{array}%
\right],
\end{equation}
or as a single (exact) equation for 1-forms,
\begin{equation}\label{form1}
    dv=\pm\sqrt{1+\frac{\nu^2}{u_x^2+u_y^2}}\left(-u_ydx+u_xdy\right).
\end{equation}
This implies the local existence of a solution to the
divergence-form equation
\begin{equation}\label{minsurf1}
    \frac{\partial}{\partial
    x}\left(\sqrt{1+\frac{\nu^2}{u_x^2+u_y^2}}u_x\right)+\frac{\partial}{\partial
    y}\left(\sqrt{1+\frac{\nu^2}{u_x^2+u_y^2}}u_y\right)=0,
\end{equation}
whenever condition (\ref{poshyp1}) is satisfied.

The equations (\ref{diff}), (\ref{orthog}), and (\ref{propo}) can
also be solved for $u_x$ and $u_y,$ in addition to being solvable
for $v_x$ and $v_y$ as in (\ref{exact1}) and (\ref{exact2}). Under
the hypothesis that either
\begin{equation}\label{idhyp}
    v_x^2+v_y^2 = \nu^2
\end{equation}
or
\begin{equation}\label{poshyp2}
    v_x^2+v_y^2>\nu^2,
\end{equation}
one obtains, by completely analogous arguments to those applied to
$v_x$ and $v_y,$ the vector equation
\begin{equation}\label{vec2}
    \left[%
\begin{array}{c}
  u_x \\
  u_y \\
\end{array}%
\right]=\mp\sqrt{1-\frac{\nu^2}{v_x^2+v_y^2}}\left[%
\begin{array}{c}
  -v_y \\
   v_x \\
\end{array}%
\right]
\end{equation}
and the equation
\begin{equation}\label{form2}
    du =
    \mp\sqrt{1-\frac{\nu^2}{v_x^2+v_y^2}}\left(-v_ydx+v_xdy\right),
\end{equation}
for 1-forms. Note that eq.\ (\ref{form2}) is exact if either
(\ref{idhyp}) holds, or if (\ref{poshyp2}) holds and the
divergence-form equation
\begin{equation}\label{gasdyn1}
    \frac{\partial}{\partial
    x}\left(\sqrt{1-\frac{\nu^2}{v_x^2+v_y^2}}v_x\right)+\frac{\partial}{\partial
    y}\left(\sqrt{1-\frac{\nu^2}{v_x^2+v_y^2}}v_y\right)=0
\end{equation}
is satisfied.

\subsubsection{Hodge-B\"acklund interpretation}

Equations (\ref{minsurf1}) and (\ref{gasdyn1}) define a B\"acklund
transformation $u\rightarrow v$ and its inverse $v\rightarrow u$
(\emph{c.f.} \cite{MT3}). These equations can be written in the form
of nonlinear Hodge equations
\begin{equation}\label{nh1}
    d*\left(\rho(Q)\omega\right)=0,
\end{equation}
\begin{equation}\label{nh2}
    d\omega = 0
\end{equation}
(that is, as eqs.\ (\ref{HF1}), (\ref{HF2}) for $\Gamma\equiv 0$).
Either
\begin{equation}\label{rhomin}
    \rho(|\omega|^2)=\sqrt{1+\frac{\nu^2}{|\omega|^2}}
\end{equation}
corresponding to eq.\ (\ref{minsurf1}), or
\begin{equation}\label{rhogas}
    \hat\rho(|\xi|^2)=\sqrt{1-\frac{\nu^2}{|\xi|^2}}
\end{equation}
corresponding to eq.\ (\ref{gasdyn1}). In either case we assume that
$Q$ does not vanish, by analogy with the equations (\ref{poshyp1}),
(\ref{vanish}), (\ref{idhyp}), and (\ref{poshyp2}).

Take $\rho$ as in (\ref{rhomin}). If the domain is simply connected,
then eq.\ (\ref{nh2}) implies, analogously to Sec.\ 5.1, that there
is a 0-form $\tau$ such that
\[
    \omega = du,
\]
and a 0-form $v$ such that
\begin{equation}\label{hbac}
    dv =
\pm*\left(\sqrt{1+\frac{\nu^2}{|\omega|^2}}\omega\right)=\pm*\left(\sqrt{|du|^2+\nu^2}\frac{\omega}{|\omega|}\right).
\end{equation}
Because the Hodge operator is an isometry $-$ which is illustrated
locally by (\ref{propo}), we have
\[
|dv|^2 = |du|^2+\nu^2.
\]
Thus the Hodge-B\"acklund transformation (\ref{hbac}) yields an
invariant form of (\ref{diff}). Unlike classical B\"acklund
transformations of the eikonal equation, in this case the
Cauchy-Riemann equations are not satisfied. Rather,
\[
u_x=\mp\rho(Q)v_y
\]
and
\[
u_y=\pm\rho(Q)v_x,
\]
which is sufficient for the orthogonality condition (\ref{orthog}).

Now take $\hat\rho$ as in (\ref{rhogas}). Arguing as before, we
conclude that there is a 0-form $\tilde u$ such that $\omega =
d\tilde u,$ and a 0-form $\tilde v$ such that
\[
d\tilde v =
\pm\ast\left(\sqrt{Q-\nu^2}\frac{\omega}{|\omega|}\right).
\]
We obtain
\[
|d\tilde v|^2 = |d\tilde u|^2-\nu^2.
\]
Letting $\tilde u = \pm i u$ and $\tilde v = \pm i v,$ we obtain a
mapping taking solutions to eqs.\ (\ref{HF2}), (\ref{HB}) with
$\Gamma\equiv 0$ and $\hat\rho$ satisfying (\ref{rhogas}) into
solutions of that system with $\Gamma\equiv 0$ and $\rho$ satisfying
(\ref{rhomin}).

These arguments extend immediately to the Hodge-Frobenius case $-$
for example, by replacing (\ref{hbac}) with the expansion
\[
e^{\hat\eta} dv =
\pm\ast\left(\sqrt{e^{2\eta}|du|^2+\nu^2}\frac{\omega}{|\omega|}\right),
\]
and squaring both sides. They also extend in a straightforward way
to gradient-recursive $k$-forms and, with some modifications, to
general $k$-forms (see the remarks following Theorem 10).

Motivated by these examples, we offer a general definition of the
\emph{Hodge-B\"acklund transformation}. It is a map taking a
solution $a$ of a nonlinear Hodge-Frobenius equation having mass
density $\rho_A$ into a solution $b$ of a nonlinear Hodge-Frobenius
equation having mass density $\rho_B$ and vice-versa, where $B$ may
equal $A$ but $b$ will not equal $a.$

\section{Appendix: Methods from elliptic theory}

In this appendix we collect the proofs of Theorem 5, 6, and 8, which
follow directly from the association of solutions to a uniformly
subelliptic operator via Lemmas 2-4.

\subsection{Proof of Theorem 5}

We require the following well known extension of de
Giorgi-Nash-Moser theory:

\bigskip

\textbf{Theorem A.1} (Morrey \cite{M}, Theorem 5.3.1). \emph{Let
$n
> 2.$ Let $U\in H^{1,2}(D)$ for each
$n$-disc $D\subset\subset\Omega,$ where $U(\mathbf{x})\geq 1$ and
define an $L^2$-function $W=U^\lambda$ for some
$\lambda\in\left[1,2\right).$ Suppose that}
\begin{equation}\label{mor1}
    \int_\Omega\left(a^{\alpha\beta}\partial_\beta W\partial_\alpha\zeta+fW\zeta\right)d\Omega\leq 0
\end{equation}
\emph{$\forall\zeta\in C_0^\infty(\Omega)$ with
$\zeta(\mathbf{x})\geq 0,$ where the coefficients $a$ and $f$ are
measurable with $f\in L^{n/2}(\Omega).$ Let the matrix $a$ satisfy
the ellipticity condition}
\[
C|\xi|^2\leq a^{\alpha\beta}(\mathbf{x})\xi_\alpha\xi_\beta
\]
\emph{for $|a(\mathbf{x})|\leq M$ at $a.e.\,\mathbf{x}\in\Omega$
and for all $\xi.$ Moreover, let the growth condition (\ref{mor})
be satisfied. If $U\in L^2(\Omega),$ then $U$ is bounded on each
$n$-disc $D\subset\subset\Omega$ and satisfies}
\[
|U(\mathbf{x})|^2\leq
Ca_0^{-n}\int_{D_{\left(R+a_0\right)}\left(\mathbf{x}_0\right)}|U(\mathbf{y})|^2d\mathbf{y},
\]
\emph{where $\mathbf{x}$ is a point of}
$D_R\left(\mathbf{x}_0\right)\subset\Omega.$

\bigskip

Slightly modified conditions will extend Morrey's result to $n=2;$
see \cite{M}, Sec.\ 5.4.

We apply Theorem A.1, taking $U(x)=H\left(Q(\mathbf{x})\right)+1$
and $\lambda = 1.$  Let $a$ be given by the matrix $\alpha$ of
Lemma 4 and let $f$ be given by (\ref{eff}). Then $f$ satisfies
the local growth condition (\ref{mor}) by smoothness. The
inequality of Lemma 3 completes the proof of Theorem 5.

\subsection{Proof of Theorem 6}

Without loss of generality, take $p$ to be the origin of
coordinates. Write the operator $\mathcal L_\omega$ in the form
\[
\mathcal L_\omega(H) = \Delta H + \nabla\cdot\left[T(H)\right]
\]
for
\[
T =
\left(-1\right)^{3n}\ast\left\{\omega\wedge\ast\left[\rho'(Q)dQ\wedge\omega\right]\right\}.
\]
Then
\[
\nabla\left[T(H)\right] =
\partial_k\left(a^{jk}\partial_jH\right)
\]
by arguments analogous to those of Lemma 4. Moreover, Lemma 4
implies that
\begin{equation}\label{ellip}
    |\nabla H|^2\leq\left(1+a^{jk}\right)\partial_j H \cdot\partial_k H\leq C|\nabla
    H|^2.
\end{equation}
If $f$ is given by (\ref{eff}), then
\begin{equation}\label{morrey}
    \int_\Omega\left\{\left[\nabla
H+T(H)\right]\cdot\nabla\zeta-CfH\zeta\right\}d\Omega\leq 0
\end{equation}
for all nonnegative test functions $\zeta\in C_0^\infty(\Omega)$
which vanish in a neighborhood of the origin. Define
\[
\zeta = \left(\eta\overline{\eta}_\nu\right)^2H,
\]
where $\eta\in C_0^\infty\left(D\right)$ for an $n$-disc $D$ of
radius $R,$ containing the origin, and itself completely contained
in the interior of $\Omega.$ Recalling that we have taken $H(0)=0$
and that (\ref{subsonic}) implies that $H'Q)>0,$ we conclude that
$H$ is a nonnegative function. Let $\overline\eta_\nu$ be given by
the sequence \cite{GS}
\[
\overline\eta_\nu(x)=\left\{\begin{array}{cc}
  0 & |x|\leq \nu^{-2} \\
  \log \left(\nu^2|x|\right)/\log\left( \nu^2R\right) & \nu^{-2}< |x|< R \\
  1 & R \leq|x|. \\
\end{array}
\right.
\]
Notice that $\overline\eta_\nu$ vanishes on a neighborhood of the
origin for any finite parameter $\nu;$ but as $\nu$ tends to
infinity, $\overline\eta_\nu$ converges pointwise to 1, whereas
$\nabla\overline\eta_\nu$ converges to zero in $L^n(D).$

Inequality (\ref{morrey}) now assumes the form
\[
0\leq\int_D\left[\left(\mathcal
L+f\right)H\right]\cdot\left(\eta\overline\eta_\nu\right)^2H \ast
1
\]
\[
=\int_D\left\{\left[\partial_k\left(1+a^{jk}\right)\partial_j
+f\right]H\right\}\left(\eta\overline\eta_\nu\right)^2H \ast 1 =
\]
\[
-\int\left(1+a^{jk}\right)\partial_j
H\cdot\partial_k\left[\left(\eta\overline\eta_\nu\right)^2H\right]\ast
1 + \int_D f H^2\left(\eta\overline\eta_\nu\right)^2 \ast 1,
\]
where
\[
-\int_D\left(1+a^{jk}\right)\partial_j
H\cdot\partial_k\left[\left(\eta\overline\eta_\nu\right)^2H\right]\ast
1 =
\]
\[
-2\int_D\left(1+a^{jk}\right)\overline\eta_\nu^2\eta\left(\partial_k\eta\right)
H\cdot\partial_j H\ast 1
\]
\[
-2\int_D\left(1+a^{jk}\right)\eta^2\overline{\eta}_\nu\left(\partial_k\overline\eta_\nu\right)
H\cdot\partial_j H \ast 1
\]
\[
-\int_D\left(1+a^{jk}\right)\left(\eta\overline\eta_\nu\right)^2\left(\partial_j
H\cdot\partial_k H\right)\ast 1
\]
\[
 \equiv
-\left(2i_1+2i_2+i_3\right).
\]
That is,
\[
i_3
=\int_D\left(1+a^{jk}\right)\left(\eta\overline\eta_\nu\right)^2\left(\partial_j
H\cdot\partial_k H\right)\ast 1
\]
\[
\leq
2\left(|i_1|+|i_2|\right)+\int_D\left(\eta\overline\eta_\nu\right)^2fH^2\ast
1.
\]
Estimating the integrals on the right individually, we have
\[
2i_1 =
2\int_D\left(1+a^{jk}\right)\overline\eta_\nu^2\eta\left(\partial_k\eta\right)
H\cdot\partial_j H\ast 1\leq
\]
\[
\varepsilon\int_D\left(\overline\eta_\nu\eta\right)^2|\nabla
H|^2\ast 1 +
\]
\[
\frac{1}{\varepsilon}\int_D\overline\eta_\nu^2|\nabla\eta|^2H^2\ast
1.
\]
\[
2i_2
=2\int_D\left(1+a^{jk}\right)\eta^2\overline{\eta}_\nu\left(\partial_k\overline\eta_\nu\right)
H\cdot\partial_j H \ast 1\leq
\]
\[
\varepsilon\int_D\left(\overline\eta_\nu\eta\right)^2|\nabla
H|^2\ast 1 +
\frac{1}{\varepsilon}\int_D\eta^2|\nabla\overline\eta_\nu|^2H^2\ast
1
\]
\[
\equiv\varepsilon i_{21}+\frac{1}{\varepsilon}i_{22},
\]
where the small constants $\varepsilon$ depend on the constant $C$
in the upper inequality of (\ref{ellip}), and
\begin{equation}\label{van}
    i_{22}=\int_D\eta^2|\nabla\overline\eta_\nu|^2H^2\ast 1\leq
C||\nabla\overline\eta_\nu||_n^2||\eta
H||_{2n/\left(n-2\right)}^2.
\end{equation}
The right-hand side of (\ref{van}) tends to zero as $\nu$ tends to
infinity. Absorbing small constants on the left, now using the
lower inequality of (\ref{ellip}), we conclude that
\begin{eqnarray}
\left(1-2\varepsilon\right)\int_D\left(\overline\eta_\nu\eta\right)^2|\nabla
H|^2\ast 1\leq\nonumber\\
C(\varepsilon)\left(\int_D\overline\eta_\nu^2|\nabla\eta|^2H^2\ast
1 +||\nabla\overline\eta_\nu||_n^2||\eta
H||_{2n/\left(n-2\right)}^2+||f||_{n/2}||\eta\overline\eta_\nu
H||_{2n/\left(n-2\right)}^2\right),\label{subtr}
\end{eqnarray}
where
\begin{eqnarray}
||\eta\overline\eta_\nu H||_{2n/\left(n-2\right)}^2\leq
C\int_D|\nabla\left(\eta\overline\eta_\nu
H\right)|^2\ast 1\leq\nonumber\\
C\left(\int_D|\nabla\left(\eta\overline\eta_\nu\right)|^2H^2\ast
1+\int_D\left(\eta\overline\eta_\nu\right)^2|\nabla H|^2\ast
1\right).\label{last_est}
\end{eqnarray}
The first integral on the extreme right-hand side of inequality
(\ref{last_est}) has essentially already been estimated in
(\ref{van}), and the second can be subtracted from the left-hand
side of (\ref{subtr}) provided that its coefficient in
(\ref{subtr}), the $L^{n/2}$-norm of $f$ over $\Omega,$ is
sufficiently small. Letting $\nu$ tend to infinity, we have in the
limit
\[
\int_D\eta^2|\nabla H|^2\ast 1\leq C\int_D|\nabla\eta|^2H^2\ast 1.
\]
We conclude that $H$ is a weak solution in a neighborhood of the
singularity. This completes the proof of Theorem 6.

\subsection{Proof of Theorem 8}

\emph{1. Outline}: The result follows from the application of Lemmas
3 and 4 and Theorem 5, above, to the proofs of Theorem 6 and
Corollary 7 of \cite{O2}. The absence of condition (\ref{nonnewt})
results in a small change of $L^p$ conditions on the solution.

\emph{2. Details}: Initially proceed as Sec.\ 7.2, but choose the
test of functions $\overline\eta_\nu$ to be a sequence of
functions possessing the following properties:
\begin{itemize}
\item[(a)] $\overline\eta_\nu\in\left[0,1\right]\forall\nu;$
\item[(b)] $\overline\eta_\nu=0$ in a neighborhood of
$\Sigma\,\forall\nu;$ \item[(c)]
$\lim_{\nu\rightarrow\infty}\overline\eta_\nu=1$ a.e., \item[(d)]
$\lim_{\nu\rightarrow\infty}||\nabla\overline\eta_\nu||_{L^{n-m-\varepsilon}}=0$
\end{itemize}
(\emph{c.f.} \cite{Se2}, Lemma 2 and p.\ 73). The function $\mathcal
F$ is given by \cite{Se1}
\[
\mathcal F(H) = \left\{\begin{array}{cc}
  H^{\gamma_2}, & 0\leq H\leq\ell, \\
  \left(1/\gamma_1\right)\left[\gamma_2\ell^{\gamma_2-\gamma_1}H^{\gamma_1}-\left(\gamma_1-\gamma_2\right)\ell^{\gamma_2}\right], & H\geq\ell. \\
\end{array}
\right.
\]
The functions $\mathcal F(H)$ and
\[
G(H)\equiv \mathcal F(H)F'(H)-\gamma_2
\]
satisfy (\cite{Se1}, p. 280; see also Sec.\ 3 of \cite{GS}):
\[
\mathcal
F(H)\leq\frac{\gamma_2}{\gamma_1}\ell^{\gamma_2-\gamma_1}H^{\gamma_1};
\]
\[
H\mathcal F'(H)\leq \gamma_2\mathcal F;
\]
\[
|G(H)|\leq \mathcal F(H)\mathcal F'(H);
\]
\[
G'(H)\geq C\mathcal F'(H)^2.
\]
Replace the test function in Sec.\ 7.2 by the test function
\[
\zeta = \left(\eta\overline\eta_\nu\right)^2G(H),
\]
where $\eta$ is defined as in Sec.\ 7.2, and substitute this value
into (\ref{morrey}). We obtain
\[
0\leq\int_D\left\{\left[\partial_k\left(1+a^{jk}\right)\partial_j+f\right]H\right\}\left(\eta\overline\eta_\nu\right)^2
G(H) \ast 1=
\]
\[
-\int_D\left(1+a^{jk}\right)\partial_j
H\cdot\partial_k\left[\left(\eta\overline\eta_\nu\right)^2G(H)\right]\ast
1+\int_D f\eta^2\overline\eta_\nu^2 H\cdot G(H) \ast 1,
\]
or
\[
2\int_D\left(1+a^{jk}\right)\left(\eta\overline\eta_\nu\right)\partial_k\left(\eta\overline\eta_\nu\right)\left(\partial_jH\right)G(H)\ast
1+\]
\[
\int_D\left(1+a^{jk}\right)\left(\eta\overline\eta_\nu\right)^2G'(H)\partial_j
H\partial_k H\ast 1 \leq\int_D
f\left(\eta\overline\eta_\nu\right)^2H\cdot G(H) \ast 1.
\]
Writing this inequality in the short-hand form
\begin{equation}\label{finala}
    I_1+I_2\leq I_3,
\end{equation}
we proceed analogously to inequalities (28)-(35) of \cite{O1},
making the following changes from the notation of \cite{O1} to our
notation: $\psi\rightarrow\overline\eta,$ $Q\rightarrow H,$
$H\rightarrow \mathcal F,$ $\Xi\rightarrow G,$ $\Phi\rightarrow
f.$ Explicitly,
\[
I_1 =
2\int_D\left(1+a^{jk}\right)\left(\eta\overline\eta_\nu\right)\partial_k\left(\eta\overline\eta_\nu\right)\cdot
G(H)\partial_j H \ast 1
\]
\[
\geq
-C\int_D\left(\eta\overline\eta_\nu\right)|\nabla\left(\eta\overline\eta_\nu\right)||\mathcal
F(H)\nabla\mathcal F(H)| \ast 1
\]
\[
\geq
-\varepsilon\int_D\left(\eta\overline\eta_\nu\right)^2|\nabla\mathcal
F(H)|^2\ast 1 -
C(\varepsilon)\int_D|\nabla\left(\eta\overline\eta\right)|^2\mathcal
F(H)^2\ast 1;
\]
\[
I_2 \geq C\int_D\left(\eta\overline\eta_\nu\right)^2|\mathcal
F'(H)|^2|\nabla H|^2 \ast 1
\]
\[
\geq C\int_D\left(\eta\overline\eta_\nu\right)^2|\nabla\mathcal
F(H)|^2 \ast 1;
\]
\[
I_3 = \int_D f\left(\eta\overline\eta_\nu\right)^2H\cdot G(H) \ast
1 \leq
\]
\[
\int_D|f|\left(\eta\overline\eta_\nu\right)^2|H\cdot\mathcal
F'(H)||\mathcal F(H)|\ast 1 \leq
\]
\begin{equation}\label{finalb}
    \gamma_2\int_D|f|\left(\eta\overline\eta_\nu\right)^2|\mathcal
F(H)|^2 \ast 1 \leq
||f||_{n/2}\left(\int_D\left\vert\eta\overline\eta_\nu\mathcal
F(H)\right\vert^{2n/\left(n-2\right)}\ast
1\right)^{\left(n-2\right)/n}.
\end{equation}
As in (\ref{last_est}), we apply the Sobolev inequality to the
right-hand side of this expression, followed by the Minkowski and
Schwartz inequalities:
\[
\left(\int_D\left\vert\eta\overline\eta_\nu\mathcal
F(H)\right\vert^{2n/\left(n-2\right)}\ast
1\right)^{\left(n-2\right)/n}\leq
C\int_D|\nabla\left[\eta\overline\eta_\nu\mathcal
F(H)\right]|^2\ast 1
\]
\[
\leq
C\left[\int_D|\left(\nabla\eta\right)\overline\eta_\nu\mathcal
{F}(H)|^2\ast 1 +
\int_D|\left(\nabla\overline\eta_\nu\right)\eta\mathcal{F}(H)|^2\ast
1 + \int_D\left(\eta\overline\eta_\nu\right)^2|\nabla\mathcal
F(H)|^2\ast 1\right]
\]
\[
\equiv I_{31}+I_{32}+I_{33}.
\]
The term $I_{33}$ can be subtracted from the left-hand-side of
inequality (\ref{finala}), as its coefficient in (\ref{finalb}),
the $L^{n/2}$-norm of $f,$ is small on small discs as a
consequence of condition (\ref{mor}).
\[
I_{32}=\int_D|\left(\nabla\overline\eta_\nu\right)\eta\mathcal{F}(H)|^2\ast
1 \leq
C\left(\gamma_1,\gamma_2,\ell\right)\int_D|\left(\nabla\overline\eta_\nu\right)^2\eta^2H^{2\gamma_1}\ast
1
\]
\[
\leq
C||\nabla\overline\eta_\nu||^2_{n-m-\varepsilon}||H^{2\gamma_1}||_\beta.
\]
Letting $\nu$ tend to infinity, the term on the right-hand side is
zero for every value of $\ell.$ Now letting $\ell$ tend to
infinity and using Fatou's inequality, we conclude that
\[
\int_D\eta^2|\nabla\left(H^{\gamma_2}\right)|^2 \ast 1\leq
C\int_D|\nabla\eta|^2 H^{2\gamma_2} \ast 1.
\]
Apply Theorem A.1, taking $U=H^{\gamma_2}.$ Then $U\in H^{1,2}(D)$
and $W=U^\lambda$ satisfies inequality (\ref{mor1}) for $\lambda =
1/\gamma_2.$ Because $\gamma_2>1/2,$ we can conclude that
$\lambda<2$ as required by Theorem A.1. We now want to check that
we can choose $\gamma_2\leq 1,$ in order to obtain $\lambda \geq
1$ as also required by Theorem A.1. Because $H$ lies in the space
$L^{2\beta\gamma_1}(D)\cap L^{2\gamma_2}(D),$ we let
$\gamma_1\beta=\gamma_2.$ Substituting the definition of $\beta,$
we find that we can choose $\gamma_2\leq 1$ for $\gamma_1>1/2$
provided $m+4<n,$ which is satisfied by hypothesis. This completes
the proof of Theorem 8.

\section{Acknowledgments}

We are grateful to Richard M. Schoen for suggesting reference
\cite{T}, Marshall Slemrod for discussing the motivating example
of Sec.\ 1.2, and Yisong Yang for discussing \cite{D}.

\end{document}